\documentclass[oneside, reqno]{amsart}

\usepackage{xypic}
\input xy
\xyoption{all}
\usepackage{epsfig}
\usepackage{amsthm}
\usepackage{amssymb}
\usepackage{amsmath}
\usepackage{amscd}
\usepackage{color}
\usepackage[T1]{fontenc}
\usepackage{stackengine}
\usepackage{upgreek}
\usepackage[font=scriptsize]{caption}
\usepackage{wrapfig}
\usepackage{esvect}

\newcommand{\bg}{\begin{equation}}
\newcommand{\ed}{\end{equation}}
\newcommand{\bga}{\begin{eqnarray}}
\newcommand{\eda}{\end{eqnarray}}
\newcommand{\pf}{\textbf{Proof:\ }}

\def\cbdu{\par{\raggedleft$\Box$\par}}

\newtheorem {Theorem}  {Theorem}

\numberwithin{Theorem}{section}

\newtheorem {Lemma}[Theorem]  {Lemma}

\theoremstyle{definition}
\newtheorem{Definition}[Theorem]{Definition}
\theoremstyle{remark}
\newtheorem{Remark}[Theorem]{\bf Remark}

\expandafter\chardef\csname pre amssym.def
at\endcsname=\the\catcode`\@ \catcode`\@=11
\def\undefine#1{\let#1\undefined}
\def\newsymbol#1#2#3#4#5{\let\next@\relax
 \ifnum#2=\@ne\let\next@\msafam@\else
 \ifnum#2=\tw@\let\next@\msbfam@\fi\fi
 \mathchardef#1="#3\next@#4#5}
\def\mathhexbox@#1#2#3{\relax
 \ifmmode\mathpalette{}{\m@th\mathchar"#1#2#3}%
 \else\leavevmode\hbox{$\m@th\mathchar"#1#2#3$}\fi}
\def\hexnumber@#1{\ifcase#1 0\or 1\or 2\or 3\or 4\or 5\or 6\or 7\or 8\or
 9\or A\or B\or C\or D\or E\or F\fi}

\font\teneufm=eufm10 \font\seveneufm=eufm7 \font\fiveeufm=eufm5
\newfam\eufmfam
\textfont\eufmfam=\teneufm \scriptfont\eufmfam=\seveneufm
\scriptscriptfont\eufmfam=\fiveeufm

\catcode`\@=\csname pre amssym.def at\endcsname

\newcounter{remark}
\newcommand{\R}{\mathbf{R}}

\renewcommand{\div}{\mbox{div}}

\def  \R   {{\mathbb R}}

\def  \12  {{\frac{1}{2}}}

\def\build#1_#2^#3{\mathrel{\mathop{\kern 0pt#1}\limits_{#2}^{#3}}}

\numberwithin{equation}{section}

\begin{document}

\title[Dyadic models]{Dyadic models for fluid equations: a survey}


\author [Alexey Cheskidov]{Alexey Cheskidov}

\address{Department of Mathematics, Statistics and Computer Science, University of Illinois at Chicago, Chicago, IL 60607, USA}
\email{acheskid@uic.edu}

\author [Mimi Dai]{Mimi Dai}

\address{Department of Mathematics, Statistics and Computer Science, University of Illinois at Chicago, Chicago, IL 60607, USA}
\email{mdai@uic.edu}

\author [Susan Friedlander]{Susan Friedlander}

\address{Department of Mathematics, University of Southern California, Los Angeles, CA 90089, USA}
\email{susanfri@usc.edu} 



\begin{abstract}

Over the centuries mathematicians have been challenged by the partial differential equations (PDEs) that describe the motion of fluids in many physical contexts. Important and beautiful results were obtained in the past one hundred years, including the groundbreaking work of Ladyzhenskaya on the Navier-Stokes equations. However crucial questions such as the existence, uniqueness and regularity of the three dimensional Navier-Stokes equations remain open. Partly because of this mathematical challenge and partly motivated by the phenomena of turbulence, insights into the full PDEs have been sought via the study of simpler approximating systems that retain some of the original nonlinear features. One such simpler system is an infinite dimensional coupled set of nonlinear ordinary differential equations referred to a dyadic model. In this survey we provide a brief overview of dyadic models and describe recent results. In particular, we discuss results for certain dyadic models in the context of existence, uniqueness and regularity of solutions. 

\bigskip


\hspace{0.02cm}CLASSIFICATION CODE: 35Q35, 76D03, 76W05.
\end{abstract}

\maketitle

\begin{center}
{\Large{Dedicated to Olga Aleksandrovna Ladyzhenskaya}}
\end{center}

\medskip

\section{Introduction}
\label{sec-intro}

\subsection{Fluid equations} 
\label{sec-overview}

\noindent

The motion of incompressible fluids is mathematically described by the set of partial differential equations 
\begin{equation}\label{nse}
\begin{split}
u_t+(u\cdot\nabla) u+\nabla p=&\ \nu \Delta u,\\
\nabla\cdot u=&\ 0.
\end{split}
\end{equation}
In the system, the unknowns are the velocity field $u=u(x,t)$ and the scalar pressure $p=p(x,t)$. 
The constant $\nu$ denotes the viscosity coefficient. In the viscous case with $\nu>0$, system (\ref{nse}) is called the Navier-Stokes equation (NSE); while in the inviscid case of $\nu=0$, (\ref{nse}) is called the Euler equation. System (\ref{nse}) is posed on a spatial domain $\Omega\subset \mathbb R^n$ with either $n=2$ or $n=3$. 

The central problem that challenges mathematicians is the well-posedness of (\ref{nse}) in the Hadamard sense: existence and uniqueness of a physically reasonable solution, and the continuous dependence of solutions on the initial state. 
Two important properties have played vital roles in the study of well-posedness of the NSE: the energy balance and scaling. Upon assuming solutions are smooth, (\ref{nse}) obeys the energy law
\begin{equation}\label{basic-en}
\frac12\frac{d}{dt}\int_{\Omega} |u(t)|^2\, dx+\nu\int_{\Omega}|\nabla u(t)|^2\, dx=0.
\end{equation} 
Thus (\ref{basic-en}) implies solutions of (\ref{nse}) have the a priori energy estimates
\begin{equation}\label{Leray-space}
u\in L^\infty(0,T; L^2)\cap L^2(0,T; H^1).
\end{equation}
The NSE (\ref{nse}) has the natural scaling: if $(u(x,t), p(x,t))$ solves (\ref{nse}) with initial data $u_0(x)$, then the pair
\[u_{\lambda}(x,t)=\lambda u(\lambda x, \lambda^{2}t), \ \ p_{\lambda}(x,t)=\lambda^2 p(\lambda x, \lambda^{2}t) \]
solves (\ref{nse}) with initial data $\lambda u_0(\lambda x)$. An invariant space under such scaling is called a critical space. The critical Sobolev space for (\ref{nse}) is $\dot H^{\frac{n}{2}-1}(\mathbb R^n)$. It is easy to see that the energy space $L^2$ is critical in 2D. Hence in view of the a priori energy estimates, the 2D NSE is said to be energy critical. This is one of the main reasons that the 2D NSE has been understood well with the contributions of many mathematicians, in particular Olga Ladyzhenskaya. One also notices that for the 3D NSE the critical Sobolev space is $\dot H^{\frac12}$ which is higher than the energy space $L^2$. For this reason, the 3D NSE is said to be supercritical. The essential obstacle in the analysis of the 3D NSE stems from the supercritical feature. 

It is not clear whether a classical solution exists or not for the 3D NSE (\ref{nse}) associated with general initial data. Leray introduced the concept of weak solutions in his pioneering work \cite{Le} in 1934. In particular, he showed the existence of global in time of Leray-Hopf solutions (see also Hopf's work \cite{Hopf}) which satisfy (\ref{nse}) in the distributional sense, belong to the spaces in (\ref{Leray-space}) and satisfy the basic energy inequality. The mathematical study of the NSE has been continued ever since. Nevertheless, some basic questions regarding the well-posedness problem remain open. We shall discuss more on this in Subsection \ref{sec-open}.



\subsection{Ladyzhenskaya's contributions in hydrodynamics}
\label{sec-Lady}

\noindent

Among many fundamental contributions made in partial differential equations, Ladyzhenskaya achieved significant results in the investigation of the incompressible NSE (c.f. \cite{Lady5}), which was her favorite topic. We highlight a few of her remarkable results on the unique solvability and attractors of the NSE in the following discussion. 

The initial boundary value problem for the 3D NSE was studied by Ladyzhenskaya and Kiselev \cite{LadyK} in 1957. They showed the existence of a unique solution of the problem with Dirichlet boundary condition in the class $L^\infty(0,T; L^4(\Omega))\cap W^{1,2}(\Omega\times [0,T])$. The maximum existence time $T$ depends on the $W^{2,2}$ norm of the initial data; in particular, $T=\infty$ if the initial data is small enough in $W^{2,2}$. In her monograph \cite{Lady5}, Ladyzhenskaya pointed out that a Leray-Hopf weak solution for the 3D NSE is unique in the class $L^8(0,T; L^4(\Omega))$. It is a particular case of the Prodi-Serrin \cite{Prodi, Serrin} condition on the uniqueness of Leray-Hopf solution,
\begin{equation}\label{Prodi}
u\in L^p(0,T; L^q(\Omega)),  \ \ \ \ \frac{2}{p}+\frac{3}{q}\leq 1, \ \ \ q>3.
\end{equation}
Ladyzhenskaya \cite{Lady1} further proved that a Leray-Hopf solution satisfying (\ref{Prodi}) is in fact smooth. Condition (\ref{Prodi}) is also called the Ladyzhenskaya-Prodi-Serrin regularity criterion. We note that the end point criterion $L^\infty(0,T; L^3(\Omega))$ was finally obtained by Escauriaza, Seregin and \v Sver\'ak \cite{ESS} in 2003.

Another marvellous result of Ladyzhenskaya is that she showed in \cite{Lady0} the existence of a unique global in time solution to the 2D NSE with appropriate initial and boundary conditions. About the same time, Lions and Prodi \cite{LP} obtained the same result by a different approach; interestingly, they used the Ladyzhenskaya inequality
\begin{equation}\notag
\|u\|^2_{L^4(\Omega)}\leq C \|\nabla u\|_{L^2(\Omega)}\|u\|_{L^2(\Omega)}
\end{equation} 
for $\Omega\subset\mathbb R^2$ and a universal constant $C>0$.

Ladyzhenskaya also made outstanding contributions to the theory of attractors for the NSE. She first studied the semigroup generated by the initial boundary value problem of the 2D NSE \cite{Lady-att1}. She provided a definition of the global attractor and showed its existence for the 2D NSE. Moreover, she discovered some important properties of the global attractor, for instance, it is finite-dimensional. In the monograph \cite{Lady-att2}, Ladyzhenskaya further investigated the problem of attractors for other dissipative systems, including parabolic ones using the semigroup approach she developed in \cite{Lady-att1}.


\subsection{Open questions remained}
\label{sec-open}

\noindent

In spite of the extensive efforts of many generations of mathematicians, including Ladyzhenskaya, important questions for the NSE still remain open. One of the notable open problems is whether a classical solution for the 3D NSE exists for all the time or it develops singularity in finite time. As mentioned above, this problem for the 2D NSE (which is critical) was affirmatively answered by Ladyzhenskaya. However, for the 3D NSE which is supercritical, the delicate competition between the nonlinear term and the linear dissipative term poses great obstacles to find the answer to the global regularity problem. In attempts to gain insights for this challenging problem, approximating models have been introduced and studied. For instance, Tao \cite{Tao} modified the 3D NSE based on a dyadic model and constructed finite time blow-up solutions for the modified NSE. 

Another unanswered question is the (non-)uniqueness of Leray-Hopf solution for the 3D NSE. In the joint paper \cite{LadyK} with Kiselev, Ladyzhenskaya pointed out that the Leray-Hopf class is too wide to prove uniqueness.  Relatively recently, significant progress was made on this topic. Jia and \v Sver\'ak \cite{JS1} provided sufficient conditions such that there are non-unique Leray-Hopf solutions for the 3D NSE. Buckmaster and Vicol \cite{BV} constructed non-unique weak solutions with finite energy for the 3D NSE, by extending convex integration techniques from the Euler equation to the NSE. Cheskidov and Luo \cite{CL} constructed non-unique weak solutions in spaces close to the border of the Prodi-Serrin regime represented by (\ref{Prodi}). With external forcing, Albritton, Bru\'e and Colombo \cite{ABC} constructed two distinct Leray-Hopf solutions with the help of instability analysis of a certain background solution. We note that the forcing approaches infinity when the time $t\to 0$. 

We denote a solution of the NSE (\ref{nse}) with $\nu>0$ by $u^\nu$ and a solution of the Euler equation (\ref{nse}) with $\nu=0$ by $u$. It is natural to ask whether $u^\nu$ converges to $u$ in certain space as $\nu\to 0$. This is the vanishing viscosity limit problem which is not completely resolved yet, despite some progress made in very special cases. A closely related question is the vanishing viscosity limit in Kolmogorov's ``zeroth law of turbulence'' referred to as anomalous dissipation. Define the average dissipation rate as
\begin{equation}\label{dissipation-e}
\varepsilon^\nu= \lim_{T\to\infty}\frac{1}{T}\int_0^T \nu \|\nabla u^\nu(t)\|^2_{L^2}\, dt.
\end{equation}
In his 1941 seminal work \cite{K41}, Kolmogorov conjectured that 
\begin{equation}\notag
\lim_{\nu\to 0} \varepsilon^\nu =\varepsilon^0>0
\end{equation}
which indicates anomalous loss of energy in the vanishing viscosity limit. Kolmogorov also predicted that there exists a critical dissipation wavenumber $\kappa_d$ which separates the so-called inertial range from the dissipation range. Statistical tools and scaling analysis suggested 
\begin{equation}\label{kappa}
\kappa_d \sim \left(\frac{\varepsilon^\nu}{\nu^3}\right)^{\frac14}
\end{equation}
and the scaling law for the energy spectrum in the inertial range
\begin{equation}\label{spect}
\mathcal E(k)\sim (\varepsilon^\nu)^{\frac23} k^{-\frac53}, \ \ \ \mbox{for} \ \ k\leq \kappa_d.
\end{equation}
It is widely believed that properties of stationary solutions for the NSE and Euler equation will provide important insights into the understanding of the ``zeroth law of turbulence''. Naturally the topic of a global attractor also falls into the constellation of these problems. 
However our present knowledge of the global attractor for the 3D NSE is incomplete.

\subsection{Toy models}

\noindent

As an attempt to understand the such challenging problems for the NSE and Euler equation, approximating models have been introduced and investigated by many mathematicians. In fact, Ladyzhenskaya was among those who first introduced a modified NSE model in 1960s, see \cite{Lady4}. Her model is obtained by modifying the classical NSE where the velocity fluctuates rapidly in the high frequency range. With reasonably mild assumptions on the initial data, Ladyzhenskaya showed global unique solvability for her model. In early 1970s, a type of discretized approximating models became popular in the physics community, known as shell models or dyadic models. They were deployed to understand energy cascade of turbulent flows, with attention on numerical simulation of these models. Around the time of 2000s, a renewed interest in dyadic models arose among mathematicians. 
Mathematical analysis produced answers to the challenging problems of global regularity, global attractors and anomalous dissipation in the context of dyadic models. In this survey we provide a review of the topic of dyadic models.\\





\bigskip

\section{Overview on dyadic models for the fluid equations}

\subsection{Dyadic models for the NSE and Euler equation}


\noindent

As approximations of the classical NSE and Euler equations, dyadic models of the NSE and Euler are infinite systems of ordinary differential equations. Although with spatial structure and geometry over simplified, the dyadic models retain certain essential features of the NSE and Euler equation. A common property of all the dyadic models is that only local interactions between the modes in Fourier space are taken into account. As a result, compared to the original PDEs, there are considerable computational and analytical advantages for treating the dyadic models.

Some early dyadic models appeared in the pioneering works of Lorenz \cite{Lorenz}, Siggia \cite{Sig}, Desnyansky and Novikov \cite{DN}, and Obukhov \cite{Ob}. 
One such early model was introduced by Gledzer \cite{Gle} in 1970s and then generalized by Ohkitani and Yamada \cite{OY} in 1980s, which is referred as the GOY shell model. This model describes the time evolution of complex-valued Fourier components $u_n$ of a scalar velocity field, taking the form
\begin{equation}\label{goy}
\frac{d}{dt} u_n+\nu k_n^2u_n=k_n B_n[u,u]+f_n
\end{equation}
where $k_n=k_0\lambda^n$ denotes the scalar wavenumber of the $n$-th shell with $\lambda$ being the intershell ratio, and $f_n$ represents the forcing on the $n$-th shell.  The quadratic nonlinear operator $B_n$ is chosen such that 
the energy $E(t)=\sum_{n} |u_n(t)|^2$ is formally conserved for the model. Another quadratic quantity without a positive sign (sometimes referred as the helicity) is also conserved. Reviews of early dyadic models in applications investigating energy cascade mechanisms for turbulent flows can be found in \cite{Bif, BJPV}. 

Although with definite advantages in numerical studies, the GOY model presents some obstacles for deep analytical investigation.  In the late 1990s
L'vov, Podivilov, Pomyalov, Procaccia and Vandembroucq \cite{LPPPV} proposed a new shell model, called the Sabra model. This model is similar to the GOY model in many aspects: complex-valued Fourier components, conservation of energy and helicity, only local interactions of the modes, but with the bilinear operator $B_n$ in (\ref{goy}) slightly modified. The resulting fundamental difference is that the Sabra model  contains fewer complex conjugation operators in the nonlinear terms and hence exhibits shorter-ranged correlations than the GOY model. In \cite{CLT}, Constantin, Levant and Titi performed an analytic study of the Sabra model. They showed global regularity of solutions, the existence of a finite dimensional global attractor and globally invariant inertial manifolds. Moreover, they obtained the existence of an exponentially decaying energy dissipation range with sufficiently smooth forcing on large scales. Nevertheless, further analytical studies are still rather difficult. One main reason is that the Sabra model exhibits both forward and backward energy cascade mechanisms. 

In early 2000s Dinaburg and Sinai \cite{DS} suggested an approximating model for the 3D NSE which is a system of quasilinear equations in Fourier space. The model has characteristics along which the nonlinearity propagates. Denote $\vec k=(\vec k_1, \vec k_2, \vec k_3)$ by the vector wavenumber in Fourier space and $\widehat u(\vec k, t)=\left(\widehat u_1(\vec k, t), \widehat u_2(\vec k, t),\widehat u_3(\vec k, t)\right)$ by the Fourier transform of the velocity field $u(x,t)$. The approximating model is given by
\begin{equation}\label{DS}
\begin{split}
&\left[\frac{d}{dt}+\left(\vec k_m\int \vec k_n'\widehat u_m(\vec k', t)\, d\vec k'\right)\frac{\partial}{\partial \vec k_n}\right] \widehat u_m\\
=&\ \widehat u_m \int \vec k_m'\widehat u_m(\vec k', t)\, d\vec k'-\frac{2\vec k}{|\vec k|^2}\left(\widehat u_m\vec k_n \int \vec k_m'\widehat u_n(\vec k', t)\, d\vec k'\right)-\nu |\vec k_m|^2\widehat u_m,\\
\vec k\cdot \widehat u=&\ 0.
\end{split}
\end{equation}
Note that the Fourier components $\widehat u$ are real-valued vector fields for this model. Under the assumption that the initial Fourier coefficients $\widehat u(\vec k, 0)$ are non-zero only for a finite set of modes, system (\ref{DS}) reduces to a finite system of ODEs. In a particular case with two non-zero modes, Dinaburg and Sinai showed existence of solutions to the system. They pointed out that there is evidence that the asymptotic behavior of the approximating system (\ref{DS}) is consistent to that of the original NSE.

A major drawback of the Dinaburg-Sinai model is that the energy balance law is not preserved. To restore the Euler equation property of energy conservation, Friedlander and Pavlovi\'c \cite{FP} modified Dinaburg-Sinai's model such that the bilinear operator of the quadratic nonlinear part has a skew symmetry property. As a result energy is formally conserved. Consider the finite set of vectors $\{ u^{j}(t)\}$ (essentially Fourier coefficients, with $\widehat \cdot$ dropped for simplifications) associated with the wavenumber vectors $\{\vec k^j(t)\}$. Denote the matrix 
\[A^{j-1}=u^{j-1} (\vec k^j)^T\]
and 
\[B^j=B^j(u^j, \vec k^j)=\frac12\left[u^j(\vec k^j)^T+\vec k^j (u^j)^T\right].\] 
The modified vector model for the Euler equation reads
\begin{equation}\label{FP}
\begin{split}
\frac{d}{dt} u^j=&\ \lambda_j^{\frac52}\left[ A^{j-1} u^{j-1}-\frac{\left((\vec k^j)^T A^{j-1} u^{j-1}\right) \vec k^j}{(\vec k^j)^2}\right]\\
&-\lambda_{j+1}^{\frac52}\left[B^{j+1} u^j-\frac{2\left( (\vec k^j)^T B^{j+1} u^{j}\right) \vec k^j}{(\vec k^j)^2}\right],\\
(\vec k^j)^T u^j=&\ 0
\end{split}
\end{equation}
for all $j\geq 0$, with $\lambda_j=\lambda^j$ for fixed $\lambda>1$.  Within the two brackets of the right hand side of (\ref{FP}), the second terms play the role of pressure and guarantee the incompressibility. For a particular set of initial data the vector model (\ref{FP}) reduces to a scalar model given by the infinite ODE system 
\begin{equation}\label{FP-s}
\frac{d}{dt} a_j=\lambda_j^{\frac52}a_{j-1}^2-\lambda_{j+1}^{\frac52}a_ja_{j+1}
\end{equation}
where the scalar-valued function $a_j(t)$ represents an average of the vector-valued function $u_j(t)$. One notices that 
the scalar model (\ref{FP-s}) has a built in forward energy cascade mechanism.  Local existence of solutions in Sobolev space $H^s$ with $s>\frac52$ was obtained in \cite{FP}. Moreover, the authors showed that finite time blow-up occurs in the sense that the $H^{\frac32+\varepsilon}$ norm becomes unbounded in finite time. More details of the analysis of (\ref{FP-s}) will be provided  later in this review. 

In the study of the partial regularity for the NSE with hyper-dissipation, Katz-Pavlovi\'c \cite{KP-cheap} introduced a dyadic model to illustrate their main ideas. The model is obtained by considering the wavelet expansion of the vector valued velocity function $u(t)$ in the spirit of Desnyansky-Novikov model and GOY model. The evolution of the wavelets coefficients satisfies a system which is essentially the same as the scalar model (\ref{FP-s}). In \cite{KP} the authors further showed that the dyadic Euler model blows up in finite time, and the dyadic NSE model with sufficiently weak dissipation blows up in finite time as well. 


To end this section we briefly mention that other dyadic models have been suggested with the aim of restoring the self-similar symmetry and intermittency feature of turbulent flows, for instance, see Mailybaev \cite{Mai0}. The author derived a nonlinear spatiotemporal scaling symmetry of inviscid hydrodynamic equations based on a simplified dyadic model. Numerical analysis was performed to confirm the hidden scale invariance of intermittent turbulence in this dyadic model. In \cite{Mai2} a shell model with complex variables was introduced and solved analytically. It was shown that the model features intermittent behavior with anomalous power-law scaling of the structure functions and its solution associates intermittency with the hidden symmetry of the turbulence. It was further illustrated in \cite{Mai1} that the intermittent dynamics is related to the self-similarity property of turbulence through the study of the Sabra dyadic model.

\medskip


\subsection{Intermittency} \label{subsec:Intermittency}

\noindent

The idea of intermittency goes back to Landau's remark questioning the exactness of Kolmogorov's laws. Intuitively, the intermittency dimension $d$ at scale $\ell$ is a number between $0$ and $3$, so
that 
\[
\text{Number of eddies} \sim \left(\frac{L}{\ell}\right)^{d}.
\]
Here $L$ is the size of the domain.
The case $d=3$ corresponds to Kolmogorov's regime
where eddies occupy the whole region for each scale in the inertial range. The other borderline 
case is $d=0$ (extreme intermittency), where the number of eddies is of order one on all the scales. Intermittent flows ($d<3$) exhibit deviations from Kolmogorov's scaling laws, which can be
measured in numerical simulations and experiments. For instance,
a direct numerical simulation performed by Kaneda et al. \cite{kaneda}  on the Earth Simulator suggests that $d \approx 2.7$ for turbulent flows.

Mathematically, we can associate intermittency with the level of saturation of Bernstein's
inequality. Let $L$ be the domain length scale. 
Denote $\lambda_j=2^j/L$ for integers $j\geq 0$. Let $u_j$ be the $j$-th Littlewood-Paley projection of a vector field $u$. 
Recall that Bernstein's inequality in three dimensional space takes the form
\begin{equation}\notag
\|u_j\|_{L^q}\leq c \lambda_j^{3(\frac1p-\frac1q)}\|u_j\|_{L^p}.
\end{equation}
Adapting the idea of \cite{CD-Kol}, we define the intermittency dimension $d$ for a 3D turbulent vector field $u$ in view of the saturation level of Bernstein's inequality,
\begin{equation}\label{int-def}
d:=\sup\left\{s\in\R:\left<\sum_{j}\lambda_q^{-1+s}\|u_j\|_{L^\infty}^2\right>\leq c^{3-s}L^{-s}\left<\sum_{j}\lambda_q^2\|u_j\|_{L^2}^2\right>\right\},
\end{equation}
where $c$ is an absolute constant.
Thus we have $d\in[0,3]$ and the optimal Bernstein's relationship 
\begin{equation}\label{bern}
\|u_j\|_{L^\infty}\sim \lambda_q^{(3-d)/2}\|u_j\|_{L^2}
\end{equation} 
at each scale $\lambda_j$. 

For instance, if the velocity of a flow on a torus $[0,L]^3$ at length scale $\ell$ looks like $u\sim(\sin(y/\ell),0,0)$, then the $L^p$ norms are of the same order, and we can argue that the eddies at the length scale $\ell$ occupy the whole region. This is called Kolmogorov's regime
and the intermittency dimension $d=3$ here. On the other hand, if $u$ is a Dirichlet kernel at
length scale $\ell$, then $\|u\|_{L^2} \sim \ell^{3/2} \|u\|_{L^\infty}$, i.e., Bernstein's inequality
is saturated, and we can see that the fluid velocity is large only in a box of length $\ell$,
i.e., there is only one eddy. This is the extreme intermittency regime and $d=0$.

The effect of intermittency on the regularity properties of solutions to the NSE (\ref{nse}) and toy models has been also studied in the past decade \cite{CF, CSreg, CSint}. Discontinuous weak solutions in the largest critical space and even supercritical spaces near $L^2$ were obtained in \cite{CD-discont, CD-norm} using  Beltrami type flows with the intermittency dimension $d=0$. Such an extreme intermittency was achieved using Dirichlet kernels.
 Roughly speaking, in order for the $n$-dimensional Navier-Stokes equations to develop singularities, the intermittency dimension $d$ of the flows should be less than $n-2$, so that  Bernstein's inequality is highly saturated. So $d=1$ is critical for the 3D NSE.

\subsection{A ``derivation'' of scalar model accounting for the intermittency effect}

\noindent

As discussed briefly above, one can obtain the scalar model (\ref{FP-s}) by different approximating approaches. We will present one way to arrive at a scalar model as (\ref{FP-s}) in which we relate the nonlinear scaling index to the intermittency effect. 

We will derive the following scalar dyadic model for the NSE
\begin{equation}\label{shell-nse}
\begin{split}
\frac{d}{dt}a_j+\nu\lambda_j^2a_j&-\alpha\left(\lambda_{j-1}^{\frac{5-d}{2}}a_{j-1}^2-\lambda_j^{\frac{5-d}{2}}a_ja_{j+1}\right)\\
&-\beta\left(\lambda_{j-1}^{\frac{5-d}{2}}a_{j-1}a_j-\lambda_{j}^{\frac{5-d}{2}}a_{j+1}^2\right)=0
\end{split}
\end{equation}
with $a_j(t)\sim \|u_j(t)\|_{L^2}$ and constants $\alpha, \beta\geq0$.
We start from the energy balance in the $j$-th shell:
\begin{equation}\label{energy-j}
\frac12\frac{d}{dt}\|u_j\|_{L^2}^2+\int_{\mathbb R^3}(u\cdot\nabla u)_j\cdot u_j\, dx+\nu\|\nabla u_j\|_{L^2}^2=0,
\end{equation}
which is obtained by projecting the NSE (\ref{nse}) onto the $j$-th shell, taking dot product with $u_j$, and integrating over the space $\mathbb R^3$ or $\mathbb T^3$. The next step is to analyze the flux through the $j$-th shell,
\begin{equation}\notag
\Pi_j := \int_{\mathbb R^3}(u\cdot\nabla u)_j\cdot u_j\, dx.
\end{equation}
We make the assumption on local interactions that only the nearest shells interact with each other. On the other hand, we notice that $\int_{\mathbb R^3}(u_i\cdot\nabla u_j)\cdot u_j\, dx=0$ for any $i,j\geq 0$ due to the divergence free property $\nabla\cdot u_j=0$. Therefore, we are able to list the non-vanishing terms in the flux as 
\begin{equation}\notag
\begin{split}
\Pi_j=&\ \int_{\mathbb R^3}(u_j\cdot\nabla u_{j+1})\cdot u_j\, dx+\int_{\mathbb R^3}(u_{j+1}\cdot\nabla u_{j+1})\cdot u_j\, dx\\
&+\int_{\mathbb R^3}(u_{j-1}\cdot\nabla u_{j-1})\cdot u_j\, dx+\int_{\mathbb R^3}(u_j\cdot\nabla u_{j-1})\cdot u_j\, dx\\
=&\ \int_{\mathbb R^3}(u_j\cdot\nabla u_{j+1})\cdot u_j\, dx+\int_{\mathbb R^3}(u_{j+1}\cdot\nabla u_{j+1})\cdot u_j\, dx\\
&-\int_{\mathbb R^3}(u_{j-1}\cdot\nabla u_{j})\cdot u_{j-1}\, dx-\int_{\mathbb R^3}(u_j\cdot\nabla u_{j})\cdot u_{j-1}\, dx
\end{split}
\end{equation}
where in the second step we applied integration by parts to the first and second integrals.
We denote 
\[Q_j=\int_{\mathbb R^3}(u_j\cdot\nabla u_{j+1})\cdot u_{j}\, dx, \ \ \ P_j=\int_{\mathbb R^3}(u_{j+1}\cdot\nabla u_{j+1})\cdot u_{j}\, dx.\]
Thus, the flux $\Pi_j$ can be rewritten as 
\[\Pi_j=Q_j-Q_{j-1}+P_j-P_{j-1}.\]
Assume $Q_j\geq0$ and $P_j\geq0$ for all $j\geq -1$. The terms $Q_j$ and $P_j$ are regarded as the energy escaping to the next shell, while $Q_{j-1}$ and $P_{j-1}$ are regarded as energy coming from the previous shell. It is important to note that in the inviscid case $\nu=0$, the total energy $\|u(t)\|_{L^2}^2=\sum _{j\geq -1}\|u_j(t)\|_{L^2}^2$ is conserved. 

Next we estimate the size of $Q_j$ and $P_j$ by using Bernstein's relation (\ref{bern}). It follows from integration by parts, H\"older's inequality, and (\ref{bern}) that
\begin{equation}\notag
\begin{split}
Q_j= \int_{\mathbb R^3}(u_j\cdot\nabla u_{j})\cdot u_{j+1}\, dx
\lesssim &\ \|u_j\|_{L^2}\|\nabla u_j\|_{L^\infty}\|u_{j+1}\|_{L^2}\\
\sim &\ \lambda_j^{\frac{5-d}{2}}\|u_j\|_{L^2}^2\|u_{j+1}\|_{L^2}
\end{split}
\end{equation}
\begin{equation}\notag
\begin{split}
P_j= \int_{\mathbb R^3}(u_{j+1}\cdot\nabla u_{j})\cdot u_{j+1}\, dx
\lesssim &\ \|\nabla u_j\|_{L^\infty}\|u_{j+1}\|_{L^2}^2\\
\sim &\ \lambda_{j}^{\frac{5-d}{2}}\|u_j\|_{L^2}\|u_{j+1}\|_{L^2}^2.
\end{split}
\end{equation}
 Define $a_j(t)\sim \|u_j(t)\|_{L^2}$. We can approximate $Q_j$ and $P_j$ as
\[Q_j=\alpha \lambda_j^{\frac{5-d}{2}}a_j^2a_{j+1}, \ \ P_j=\beta \lambda_j^{\frac{5-d}{2}}a_ja_{j+1}^2\]
for some constants $\alpha\geq 0$ and $\beta\geq 0$.
 Motivated by (\ref{energy-j}) and the analysis above, we consider the approximating equation
\begin{equation}\notag
\frac12\frac{d}{dt}a_j^2+Q_j-Q_{j-1}+P_j-P_{j-1} +\nu\lambda_j^2a_j^2=0,
\end{equation}
which leads to the shell model (\ref{shell-nse}). For the total energy $a(t)^2=\sum_{j\geq0}a_j^2$ of the approximating model, we also have the energy law
\[\frac12\frac{d}{dt}a^2+\nu \sum_{j\geq-1}\lambda_j^2a_j^2=0\]
which indicates energy conservation for smooth solutions in the inviscid case $\nu=0$.

\begin{Remark}\label{rk1}
\noindent
\begin{itemize}
\item [(i)] For 3D flows we have the intermittency dimension $d\in[0,3]$. Thus the nonlinear scaling index $\frac{5-d}{2}$ in (\ref{shell-nse}) belongs to $[1, \frac52]$ which is considered to be the physically relevant regime.
\item [(ii)] We notice that the dyadic model (\ref{shell-nse}) with $\alpha=1$, $\beta=0$ and $=0$ corresponds to the model (\ref{FP-s}) (also called the Desnyansky-Novikov or DN model), while (\ref{shell-nse}) with $\alpha=0$, $\beta=1$ and $d=0$ is Obukhov's model \cite{Ob}. 
\item [(iii)] Denote $\theta=\frac{5-d}{2}$. In the case of $\alpha=1$ and $\beta=0$, we have the integral form of (\ref{shell-nse})
\begin{equation}\notag
\begin{split}
a_j(t)=&\ a_j(0)e^{-\int_0^t \nu\lambda_j^2+\lambda_j^{\theta} a_{j+1}(\tau) \, d\tau}\\
&+\int_0^t e^{-\int_\tau^t \nu\lambda_j^2+\lambda_j^{\theta} a_{j+1}(s) \, ds} \lambda_{j-1}^\theta a_{j-1}^2(\tau)\, d\tau.
\end{split}
\end{equation}
It is clear from the integral form that $a_j(t)>0$ if $a_j(0)>0$ for all $j\geq 0$. This is usually called the positivity property. On the other hand, in the case of $\alpha=0$ and $\beta=-1$, the integral form of (\ref{shell-nse}) 
\begin{equation}\notag
\begin{split}
a_j(t)=&\ a_j(0)e^{-\int_0^t \nu\lambda_j^2-\lambda_j^{\theta} a_{j-1}(\tau) \, d\tau}\\
&-\int_0^t e^{-\int_\tau^t \nu\lambda_j^2-\lambda_j^{\theta} a_{j-1}(s) \, ds} \lambda_{j+1}^\theta a_{j+1}^2(\tau)\, d\tau
\end{split}
\end{equation}
immediately implies that $a_j(t)<0$ if $a_j(0)<0$ for all $j\geq 0$. Equivalently, the model has positivity for $\{-a_j\}$. 
\item [(iv)] Associated with the positivity of the DN model is the forward energy cascade mechanism, illustrated by
the fact
\begin{equation}\notag
\frac{d}{dt}\sum_{j=0}^Na_j^2(t)=-\nu \sum_{j=0}^N\lambda_j^2 a_j^2-\lambda_N^{\theta}a_N^2a_{N+1}<0 \ \ \mbox{for} \ \ a_{N+1}>0.
\end{equation}
Similarly, the positivity of $\{-a_j\}$ for the Obukhov model corresponds to a backward energy cascade mechanism, since
\begin{equation}\notag
\frac{d}{dt}\sum_{j=0}^Na_j^2(t)=-\nu \sum_{j=0}^N\lambda_j^2 a_j^2-\lambda_N^{\theta}a_Na^2_{N+1} >0 
\end{equation}
for $-a_{N}>0$ and $\nu=0$.
The positivity and unidirectional energy cascade mechanism have played vital roles in the analysis of the DN model and Obukhov model separately. Nevertheless, when both $\alpha>0$ and $\beta>0$, the model (\ref{shell-nse}) is a combination of the DN and Obukhov types. In this case the positivity is not preserved any more. Hence energy transfers to both high and low modes. Obviously, this is the most challenging case compared to either the DN model or Obukhov model. 
\end{itemize}
\end{Remark}

\bigskip

\section{Results for dyadic Euler and NSE models in the deterministic setting}

\medskip


\subsection{Dyadic Euler models}
\subsubsection{Regularity and blow-up problems}

\noindent

As we discussed in the Introduction,
 whether or not a solution of the 3D Euler equation develops singularity in finite time is an outstanding open problem. 
Nevertheless, for dyadic models of the Euler equation, it has been shown that finite time blow-up occurs in certain contexts. We will highlight some blow-up results in the following.

Both Friedlander \& Pavlovi\'c \cite{FP} and Katz \& Pavlovi\'c \cite{KP} proved finite time blow-up for the dyadic Euler model (\ref{FP-s}). 

\begin{Theorem}\label{thm-KP-blow} 
Let $\{a_j(t)\}$ be a solution to (\ref{FP-s}) with positive coefficients $\{a_j(0)\}$ at the initial time and 
\[a_j^2(0)>\lambda_j^{-3-\varepsilon}\]
for some $j>j_0$ with sufficiently large $j_0$. If $V(x,t)=\sum_{j} a_j(t) e^{i 2^j(\xi^j)^T x}$ is in $H^s$ for some $s>\frac52$, then the $H^{\frac32+\varepsilon}$ norm of $V(x,t)$ becomes unbounded in finite time.
\end{Theorem}

Here the $H^s$ norm is defined as 
\[ \|a(t)\|^2_{H^s}:= \sum_{j=0}^\infty \lambda_j^{2s} a_j^2(t).\]
The crucial ingredient of the proof is that solutions of (\ref{FP-s}) with positive initial data remain positive for all the time. As a result, one can see the energy flows up the shells from the nonlinear structure of the model. Rigorously, the author showed that the energy of the Carleson box 
\[E_{\mathcal B(j)}=\sum_{l\geq j} a_l^2(t)\]
monotonically increases in time. Moreover, the authors quantified the increase of the energy within a short time which allowed an iterative process to show energy concentration at high modes. Indeed, such positivity and monotonicity would be used in most of the later works. However we emphasize that such positivity is not a feature of the original Euler equation.

Later on, Kiselev and Zlato\v s \cite{KZ} showed finite time blow-up for the inviscid model 
\begin{equation}\label{model-scale1}
\frac{d}{dt} a_j=\lambda_j a_{j-1}^2-\lambda_{j+1} a_ja_{j+1}, \ \ \ j\geq 0, \ \ \lambda_j=\lambda^j, \ \lambda>1 
\end{equation}
with any non-zero initial data in $H^1$. Namely, 

\begin{Theorem}\label{thm-KZ1} 
The model (\ref{model-scale1}) with any non-zero initial data belonging to $H^1$ leads to a finite time blow-up in $H^1$. 
\end{Theorem}
We point out that the blow-up space $H^{\frac52}$ appeared in Theorem \ref{thm-KP-blow} is equivalent to the $H^1$ space of Theorem \ref{thm-KZ1}.  Indeed, taking $\lambda=2^{\frac52}$ in (\ref{model-scale1}), models (\ref{FP-s}) and (\ref{model-scale1}) are the same. The improvement from the result in Theorem \ref{thm-KP-blow} is that finite time blow-up can occur for any non-zero initial data. Thus in the inviscid case, positivity is not essential for blow-up; instead, the absence of the dissipation is the main cause of blow-up. 

The authors of \cite{KZ} also studied the Obukhov model
\begin{equation}\label{model-ob}
\frac{d}{dt} a_j=\lambda_j a_{j-1}a_j-\lambda_{j+1} a_{j+1}^2, \ \ \ j\geq 0
\end{equation}
which is equipped with a backward energy cascade mechanism. In contrast with (\ref{model-scale1}), the Obukhov model (\ref{model-ob}) has a global regular solution in $H^s$ for $s>1$. Namely,

\begin{Theorem}\label{thm-KZ2} 
For any $a^0\in H^s$ with $s>1$ and for any $T>0$, there exists a unique solution $\{a_j\}\in C([0,T], H^s)$ to (\ref{model-ob}) such that $a_j(0)=(a^0)_j$ for all $j\geq 0$.
\end{Theorem}

The authors showed an even stronger regularity result for the Obukhov model. 

\begin{Theorem}\label{thm-KZ3} 
Let $b_j(\omega)$ be independent uniformly bounded random variables such that the probability of $b_j(\omega)$ being nonpositive is uniformly bounded away from zero: 
\[P[b_j(\omega)\leq 0]>\rho>0.\]
Assume $c_j>0$ are such that $\sum_j\lambda_j^{2s}|c_j|^2<\infty$ with $s>1$. Then with probability one a solution $\{a_j(t)\}$ of the Obukhov model (\ref{model-ob}) corresponding to the initial data $a_j(0)=c_j b_j(\omega)$ satisfies $\|a\|_{H^r}\leq C(r, \omega)$ for all time $t$ and any $r<s$. Moreover, we have
\begin{equation}\notag
\lim_{t\to\infty}\|a(t)\|^2_{H^r}=\lim_{t\to\infty} a_0^2(t) = \sum_{j\geq 0}a_j^2(0).
\end{equation}
\end{Theorem}
The last statement of Theorem \ref{thm-KZ3} indicates that the solution $a(t)$ converges in $H^r$ to a constant solution with all energy concentrated in the lowest mode. This is exactly due to the backward energy cascade feature discussed in Remark \ref{rk1} [(iv)]. 

In view of Theorem \ref{thm-KZ1} and Theorem \ref{thm-KZ2}, a natural question is to investigate the global regularity problem (vice versa the blow-up scenario) for the general model 
\begin{equation}\label{model-comb}
\frac{d}{dt} a_j=\alpha\left(\lambda_j a_{j-1}^2-\lambda_{j+1} a_ja_{j+1}\right) +\beta\left(\lambda_j a_{j-1}a_j-\lambda_{j+1} a_{j+1}^2\right)
\end{equation}
with $\alpha>0$ and $\beta>0$, which is a combination of the DN type (\ref{model-scale1}) and Obukhov type (\ref{model-ob}). As we discussed in Remark \ref{rk1} [(iv)], model (\ref{model-comb}) does not preserve positivity; and the energy transfer among shells is more complicated than that of (\ref{model-scale1}) and (\ref{model-ob}) separately. Jeong and Li addressed this question in \cite{JL}. They showed that when $\beta>0$ is sufficiently small, solutions of (\ref{model-comb}) 
with nonnegative initial data can blow up in finite time. Precisely, 

\begin{Theorem}\label{thm-JL}
Let $\alpha=1$, $s>\frac13$ and $0\leq \beta< (2^s-2^{1-2s})/(1+2^{1-3s})$ in (\ref{model-comb}). Then for every initial data satisfying $a_j(0)\geq 0$ for all $j\geq 0$, the norm $\|a(t)\|_{H^s}^2$ of the solution $a(t)$ is not locally integrable on $[0,\infty)$. 
\end{Theorem} 
Note that when $s$ approaches $\frac13$, the upper bound on $\beta$ approaches 0. It is reasonable to believe that when $\beta$ is small enough, the DN type nonlinearity in (\ref{model-comb}) plays a dominant role; hence the overall energy moves up the scales and eventually concentrates on high modes to produce blow-up. We note that the authors also proved finite time blow-up for the combined model with relatively weak dissipation and more stringent condition on $\beta$. 


\medskip

\subsubsection{Anomalous dissipation, fixed point, and global attractor}
\label{subsec-fixed}

\noindent

In the statistical theory of turbulence \cite{On}, Onsager asserted that weak solutions of the Euler equation are not regular enough to conserve energy. He conjectured that the energy is conserved if the velocity field has H\"older regularity exponent $h>\frac13$ and not conserved if $h\leq \frac13$.  This is essentially connected to Kolmogorov's zeroth law -- dissipation anomaly -- discussed in Subsection \ref{sec-attrator}. The positive (conservation) side of Onsager's conjecture has been answered completely thanks to the contributions of Cheskidov-Constantin-Friedlander-Shvydkoy \cite{CCFS}, Constantin-E-Titi \cite{CET}, Eyink \cite{Ey} and Eyink-Sreenivasan \cite{ES}.  The negative (non-conservation) side of the conjecture is more challenging and was finally resolved by Isett \cite{Is} built on the previous contributions of De Lellis and Sz\'ekelyhidi \cite{DLS1, DLS2, DLS3} and collaborators, by employing the convex integration schemes. However, when the negative side of the conjecture was still open before 2007, results concerning Onsager's conjecture were obtained in the context of dyadic models.

Cheskidov, Friedlander and Pavlovi\'c \cite{CFP1, CFP2} studied the dyadic Euler model with an external forcing on large scales,
\begin{equation}\label{euler-force} 
\frac{d}{dt} a_j=\lambda_j^\theta a_{j-1}^2-\lambda_{j+1}^{\theta} a_ja_{j+1}+f_j, \ \ \ j\geq 0
\end{equation}
with $\theta=\frac{5-d}{2}$ and $d\in[0,3]$ being the intermittency dimension. 
Without loss of generality, it is assumed that 
\[f_0>0; \ \ f_j\equiv 0 \ \ \forall \ \ j\geq 1.\] 
Define
\[
\|a\|_{H^s}:=\left(\sum_j \lambda_j^{2s} a_j^2 \right)^{\frac12}, \qquad \|a\|_{B^s_{3,\infty}}:=\sup_j\lambda_j^{s+\frac{3-d}{6}} a_j,
\]
consistent with Bernstein's relationship \eqref{bern}.
A regular solution of (\ref{euler-force}) is defined to be a solution with bounded $H^{\frac13\theta}$ norm. 
The authors proved the existence of a unique fixed point for the Euler dyadic model which attracts finite energy solutions. The fixed point exhibits Onsager's scaling and does not satisfy the energy equality. 
They established the following results.


\begin{Theorem}\label{thm-CFP}
Consider the model (\ref{euler-force}) with external forcing. Then
\begin{itemize}
\item [(i)] Every regular solution satisfies the energy equality.
\item [(ii)] There exists a unique fixed point $a^*$ to (\ref{nse-force}). The fixed point is a global strong attractor and
\[a^*\notin H^{\frac13\theta}, \ \ \ a^*\in B^{\frac13}_{3,\infty} . \]
\item [(iii)] There is anomalous dissipation of energy with the rate $\varepsilon=a^*_0 f_0$. The energy spectrum of the fixed point satisfies
\[\mathcal E(k)\sim \varepsilon^{\frac23} k^{-\frac{3+2\theta}{3}}=\varepsilon^{\frac23} k^{-\frac{8-d}{3}}.\]
\item [(iv)] Every solution blows up in finite time in $H^{\frac13\theta}$ and $B^{\frac13 + \eta}_{3,\infty} $ with any $\eta>0$, and dissipates energy eventually. For $s<\frac13\theta$, the $H^s$ norm of every solution is locally square integrable. 
\end{itemize}
\end{Theorem}

The result of the fixed point being a global attractor is rather surprising considering that there is no dissipation in the system. Since the fixed point has Onsager's critical scaling but does not satisfy the energy equality, it justifies the negative side of Onsager's conjecture. In Theorem \ref{thm-CFP}, taking $d=3$, one recovers Kolmogorov's classical scaling law of the energy spectrum. As we remarked in Subsection~\ref{subsec:Intermittency}, turbulent flows exhibit intermittency, reflected in the correction to Kolmogorov's scaling law when $d<3$. The proof of the theorem relies heavily on the positivity of solutions and the monotonicity of the energy flux.

\medskip

\subsubsection{Anomalous dissipation in the unforced case}

\noindent

It is remarkable that the energy dissipation of positive solutions can be obtained for the unforced dyadic model (\ref{model-scale1}) as well. Barbato, Flandoli and Morandin proved the following results in \cite{BFM}:

\begin{Theorem}\label{thm-BFM1-posi}
Assume $a^0$ is componentwise positive and bounded. Let $a(t)$ a positive solution to (\ref{model-scale1}) with the initial data $a^0$. Then the energy of $a(t)$ 
\[E(t)=\sum_{j=0}^\infty a_j^2(t)\]
remains bounded for all positive time. 
\end{Theorem}

\begin{Theorem}\label{thm-BFM2-posi}
Let $a(t)$ be a positive solution to (\ref{model-scale1}) with finite energy. Then the energy of $a(t)$ vanishes asymptotically, i.e.
\[\lim_{t\to\infty}\sum_{j=0}^\infty a_j^2(t)=0.\]
\end{Theorem}

Theorem \ref{thm-BFM1-posi} implies that a positive solution with infinite energy initially evolves to have finite energy immediately, while Theorem \ref{thm-BFM2-posi} indicates that a positive solution with finite energy dissipates energy to zero eventually. Both of the phenomena are evidence of the anomalous dissipation for the unforced model (\ref{model-scale1}).

The authors \cite{BFM} also showed that there are initial data such that at least two solutions in $l^2$ with continuous components exist on some time interval $[0,T]$. One such solution is obtained via energy method which satisfies energy inequality. The other solution arises from a self-similar profile which increases in time.

An intriguing result on anomalous dissipation was established by Mattingly, Suidan and Vanden-Eijnden \cite{MSV} for a linear dyadic model
\begin{equation}\label{linear-model}
\frac{d}{dt} a_j=c_{j-1} a_{j-1} -c_j a_{j+1}+f_j, \ \ j\geq 0, \ \ a_{-1}\equiv 0
\end{equation}
with the external forcing $f_j$ and coupling coefficients $c_j$. The coefficients $c_j$ are chosen to ensure the energy conservation formally.  For specific $c_j$, the dynamical system (\ref{linear-model}) was shown to exhibit anomalous energy dissipation. This model illustrates some heuristic ideas similar to the random passive scalar and random Burgers equation.

\bigskip
\subsubsection{Regularizing effect of the forward energy cascade} \label{sebs:regularizing}

The basic principle proposed by Kolmogorov \cite{K41} behind turbulence is forward energy cascade. Simply put, the theory asserts that energy moves from large to small scales until it reaches the dissipation range, where the viscous forces dominate. For the Navier-Stokes equations, the dissipation range is the only tool currently used to prove regularity of solutions. However, the forward energy cascade generated by the nonlinear term might also be a mechanism to regularize solutions. For quasilinear scalar equations, the regularizing property of the nonlinear term has been studied by Tadmor and Tao in \cite{TT}, but such results remain out of reach for the Euler or Navier-Stokes equations. Tao \cite{Tao} proved blow-up for averaged Navier-Stokes equations by reducing the equations to a more complicated dyadic model where he introduced a delay in energy cascade, see Subsection~\ref{subs:AveragedNSE}. This delay is designed to destroy the regularizing effect of the nonlinear term and produces a blow-up.

Kolmogorov predicted energy cascade produces the so called dissipation anomaly, which is
characterized by the persistence of non-vanishing energy dissipation in the limit of vanishing viscosity (see Subsection~\ref{sec-attrator}). This phenomenon is related to anomalous dissipation, namely the failure of the energy to be conserved despite the absence of viscosity (see Subsection~\ref{sec-open}). Onsager conjectured that sufficiently rough solutions to Euler’s equation can exhibit anomalous dissipation, however if the solution is smooth enough, then the energy should be conserved \cite{On}. Remarkably, the energy cascade producing anomalous dissipation also leads to a gain of regularity, at least in the context of the inviscid dyadic model.
Indeed, in \cite{CFP2} Cheskidov, Friedlander, and Pavlovi\'c studied the forced inviscid dyadic model \eqref{euler-force}
showed that all the solutions of the forced inviscid dyadic model \eqref{euler-force} must have Onsager’s regularity $B^{\frac13}_{3,\infty}$ almost everywhere in time 
and  confirmed anomalous dissipation for all solutions with non-negative $l^2$ initial data
 The authors also showed that all solutions blow up in finite time in $B^{\frac13 +}_{3,\infty}$ and hence  $H^{\frac13 \theta}$ (or $H^{\frac56}$ in the case of extreme intermittency $d=3$). On the other hand, all solutions with $\ell^2$ initial data are in $B^{\frac13-}_{3,\infty}$ for almost all positive time, see Theorem~\ref{thm-CFP}. Hence, the nonlinear term regularizes solutions instantaneously.  In \cite{BM}, Barbato and Morandin studied the unforced inviscid model \eqref{euler-force} with $f=0$ and showed Onsager regularity $B^{\frac13 -}_{3,\infty}$ almost everywhere. In addition, they demonstrated that solutions remain in $B^{0}_{3,\infty}$ and hence $H^{\frac12 -\frac{d}{6}}$ for all positive time. In \cite{CZ} Cheskidov and Zaya improved  this result in the case of extreme intermittency $d=0$ ($\theta = \frac52$) by showing that  $B^{\frac{1}{10}}_{3,\infty}$ and hence $H^{\frac{3}{5}-}$ regularity is retained, which is closer to Onsager’s $B^{\frac13}_{3,\infty}$. More precisely,
\begin{Theorem}
For any positive solution to \eqref{euler-force} with $\theta = \frac52$ and initial data $a(0) \in l^2$,
\[
\sup_j \lambda_j^{\frac35}a_j(t)< \infty,
\]
for all $t > 0$.
\end{Theorem}

In fact, we conjecture that every solution must have exactly Onsager’s regularity starting some positive time, no matter whether the initial data is smooth or rough.  Hence every solution is in $B^{\frac13-}_{3,\infty}$ but not in $B^{\frac13+}_{3,\infty}$ for all $t>T$. Here $T=0$ for rough initial data $u(0) \notin B^{\frac13}_{3,\infty}$, and $T>0$ is the blow-up time for smooth (smoother than Onsager's)  initial data  $u(0) \in B^{\frac{1}{3}+}_{3,\infty}$.
\bigskip

\subsection{Dyadic NSE models}

\noindent

Now we turn to the dyadic model for the NSE 
\begin{equation}\label{nse-model1} 
\frac{d}{dt} a_j+\nu \lambda_j^{2\gamma} a_j=\lambda_j a_{j-1}^2-\lambda_{j+1} a_ja_{j+1}+f_j, \ \ \ j\geq 0
\end{equation}
with $\nu>0$ and $\gamma>0$, and $f_j$ being the external forcing for $j\geq 0$.
By rescaling the frequency we note (\ref{nse-model1}) is equivalent to the model
\begin{equation}\label{nse-force} 
\frac{d}{dt} a_j+\nu \lambda_j^{2} a_j=\lambda_j^\theta a_{j-1}^2-\lambda_{j+1}^{\theta} a_ja_{j+1}+f_j, \ \ \ j\geq 0
\end{equation}
with $\theta=\frac{1}{\gamma}$. As before, in the case of 3D scaling, we have the relation $\theta=\frac{5-d}{2}$ with $d\in[0,3]$ being the intermittency dimension. 
Both forms (\ref{nse-model1}) and (\ref{nse-force}) have been studied in the literature.

\subsubsection{Regularity and blow-up problems} \label{subs:RegBlow}

\noindent

It was first shown in \cite{KP} that solutions of (\ref{nse-model1}) without external forcing starting from positive initial data develop finite time blow-up when the dissipation degree $\gamma<\frac14$. 
Later, Cheskidov \cite{Ch} proved finite time blow-up for (\ref{nse-model1}) with $\gamma<\frac13$ if the initial data is nonnegative and large. Namely, 
\begin{Theorem}\label{thm-Ch1} 
Let $a(t)$ be a solution to (\ref{nse-model1}) with $a_j(0)\geq 0$ and $f_j=0$ for all $j\geq 0$, and $\gamma<\frac13$. Then for every $s>0$ there exists a constant $M(s)$ such that $\|a(t)\|^3_{H^{\frac13+s}}$ is not locally integrable on $[0,\infty)$ provided $\|a(0)\|_{H^{s}}>M(s)$. 
\end{Theorem}

The author also obtained global regularity for (\ref{nse-model1}) with $\gamma\geq \frac12$. 
\begin{Theorem}\label{thm-Ch2} 
Let $\gamma\geq \frac12$. For any $a(0)\in H^\gamma$, there exists a strong solution $a(t)$ to (\ref{nse-model1}) with $f_j=0$ and the initial data $a(0)$ on $[0,\infty)$.
\end{Theorem}
First note that the blow-up occurs when $\theta = \frac{1}{\gamma} >3$ or $d<-1$ in the case of 3D scaling, i.e. outside of the physically relevant range of the intermittency dimension $d\in[0,3]$. On the other hand, in the $n$-dimensional case,
\[
\theta = 1+ \frac{n-d}{2},
\]
and hence, for any $n>4$ we have that for $d>0$ small enough a blow-up occurs, as $\theta >3$ in that case. In particular, a blow-up occurs when the dyadic model scales as the five-dimensional NSE. Remarkably, 5D is a common threshold for models exhibiting a blow-up, see \cite{Tao-blog}. In fact, we conjecture that every solution of the inviscid dyadic model wants to exhibit Onsager's scaling, and solutions of the viscous dyadic model cannot have a blow-up worse than Onsager's, see Subsection~\ref{sebs:regularizing} for more discussion. Hence, the blow-up occurs only for  viscous dyadic models with critical spaces higher than Onsager's. This suggests that Theorem~\ref{thm-Ch2} is optimal.

We also note the presence of the gap $\frac13 \leq \gamma<\frac12$ between the two scenarios of finite time blow-up and global regularity. 
The gap was made smaller thanks to a regularity result of Barbato, Morandin and Romito \cite{BMR}. The authors showed global regularity for (\ref{nse-model1}) with $\gamma \geq \frac25$ by finding an invariant region for the vector $(a_j, a_{j+1})$ through a dynamical system argument. It is remarkable that $\gamma \geq \frac25$ corresponds to $\theta \leq \frac52$ and hence $d\geq0$, i.e., this covers the whole physically relevant intermittency regime. Therefore, \cite{BMR} settles  that solutions to the dyadic model corresponding to the 3D NSE are globally regular.

\medskip

\subsubsection{Global regularity in logarithmically supercritical regime}
\label{subsec-log}

\noindent

In \cite{Tao2} Tao studied the following logarithmically supercritical hyperdissipative NSE
\begin{equation}\label{nse-Tao}
\begin{split}
u_t+(u\cdot\nabla) u+\nabla p=&- D^2 u,\\
\nabla\cdot u=&\ 0,
\end{split}
\end{equation}
where $D$ is a Fourier multiplier with symbol $m:\mathbb R^n\to\mathbb R^+$.  Let $g: \mathbb R^+\to\mathbb R^+$ be a nondecreasing function such that 
\begin{equation}\label{function-g}
\int_1^\infty \frac{1}{sg^4(s)} \, ds=\infty.
\end{equation}
It is assumed that $m$ satisfies the lower bound
\begin{equation}\label{function-m}
m(\xi)\geq \frac{|\xi|^{\frac{n+2}{4}}}{g(|\xi|)}
\end{equation}
for large enough $|\xi|$. In particular, an example of such multiplier is given by
\begin{equation}\notag
m(\xi)=\frac{|\xi|^{\frac{n+2}{4}}}{\log^{\frac14}(2+|\xi|^2)}
\end{equation}
which associates with 
\begin{equation}\label{operator-D}
|D|^2=\frac{(-\Delta)^{\frac{n+2}{4}}}{\log^{\frac12}(2-\Delta)}.
\end{equation}
Note that for the operator $D$ satisfying (\ref{operator-D}), system (\ref{nse-Tao}) is slightly supercritical. 

The author showed global regularity for (\ref{nse-Tao}) with the symbol of $D$ satisfying (\ref{function-g}) and (\ref{function-m}).

\begin{Theorem}\label{thm-Tao-log} 
Let $u_0$ be smooth and compactly supported. There exists a global smooth solution to (\ref{nse-Tao}) with initial data $u_0$.
\end{Theorem}
It was conjectured by the author that the global regularity can be obtained for (\ref{nse-Tao}) with (\ref{function-g}) replaced by 
\begin{equation}\label{function-g2}
\int_1^\infty \frac{1}{sg^2(s)} \, ds=\infty.
\end{equation}
The first attempt to verify this conjecture was successfully made in the context of a dyadic model for (\ref{nse-Tao}) by Barbato, Morandin and Romito \cite{BMR3}. The authors of \cite{BMR3} proposed the following dyadic model for (\ref{nse-Tao})
\begin{equation}\label{model-log} 
\frac{d}{dt} a_j+\frac{\lambda_j}{g_j} a_j=\lambda_j a_{j-1}^2-\lambda_{j+1} a_ja_{j+1}, \ \ \ j\geq 0
\end{equation}
with $g_j$ and $\frac{\lambda_j}{g_j}$ being non-decreasing functions for all $j\geq 0$. The assumption $\sum_{j\geq 0}g_j^{-2}=\infty$ corresponds to condition (\ref{function-g}), while $\sum_{j\geq 0}g_j^{-1}=\infty$ corresponds to (\ref{function-g2}). It was shown in \cite{BMR3} that

\begin{Theorem}\label{thm-BMR-log1} 
Let $a^0\in H^s$ for all $s>0$. Assume  $\sum_{j\geq 0}g_j^{-1}=\infty$. Then there exists a solution $a(t)$ to (\ref{model-log}) with initial data $u^0$ such that $a\in L^\infty([0,\infty); H^s)$ for all $s>0$.
\end{Theorem}
The crucial element in the proof was an iterative estimate on the amount of energy stored in all shells larger than $j$ at time $t$ and the energy dissipated by these shells up to time $t$.

With insights abtained from their result for the dyadic model, the same authors provided an affirmative answer to Tao's conjecture for the original hyperdissipative NSE (\ref{nse-Tao})  in \cite{BMR2}. In particular, they showed

\begin{Theorem}\label{thm-BMR-log2} 
Let $n\geq 2$ and $u_0$ be smooth and compactly supported. There exists a global smooth solution to (\ref{nse-Tao}) associated with (\ref{function-m}) and (\ref{function-g2}) and the initial data $u_0$.
\end{Theorem}
The strategy of the proof arose from the authors' previous work \cite{BMR3} for the dyadic model of (\ref{nse-Tao}). Indeed, the smoothness of the PDE problem was reduced to the smoothness of a suitable shell model. Roughly speaking, the shell model was obtained by averaging the energy through dyadic shells. Analogously, the authors then established a recursive estimate for the contribution of energy and dissipation over large shells, which was the key ingredient to close the argument. 


\medskip

\subsubsection{Uniqueness and non-uniqueness}
\label{subsec-unique}

\noindent

It is known that a solution to the 3D NSE with high enough regularity satisfies the energy equality and is unique. In contrast energy conservation for the Euler equation, as indicated by Onsager's conjecture, the regularity threshold $\frac13$ separates the two opponent scenarios of energy conservation and anomalous dissipation. Such a threshold regularity value is not known to divide the uniqueness from non-uniqueness of solutions to the 3D NSE. 
However in the class of Leray-Hopf solutions, Ladyzhenskaya's condition $L^8(0,T; L^4(\Omega))$ and in general Prodi-Serrin's condition (\ref{Prodi}) ensure uniqueness. 

In the context of dyadic models, more knowledge on the question of uniqueness and non-uniqueness is available thanks to the work of Filonov \cite{Fi} and Filonov and Khodunov \cite{FK}. In particular, there is a critical value  of the intermittency (intrinsically connected to the regularity of the solution) that separates uniqueness from non-uniqueness scenarios. These results are described below.

We consider the dyadic model in the form (\ref{nse-force}) with the forcing $f_j$ to be specified in each case. Recall the nonlinear scaling index $\theta=\frac{5-d}{2}$ with $d\in[0,3]$ being the intermittency dimension. We define weak solution and the Leray-Hopf solution for the dyadic model in analogy with that of the original PDE.
\begin{Definition}\label{def-weak}
A sequence of functions $\{a_j(t)\}_{j=0}^\infty$ is said to be a weak solution of (\ref{nse-force}) on $[0,T)$ if $a_j\in C^1([0,T))$ satisfies (\ref{nse-force}) for all $j\geq0$, and  
\[\sum_{j=0}^\infty a_j^2(t) <\infty \ \ \ \forall \ \ t\geq0.\]
\end{Definition}

\begin{Definition}\label{def-LH}
A sequence of functions $\{a_j(t)\}_{j=0}^\infty$ is called a Leray-Hopf solution of (\ref{nse-force}) on $[0,T)$ if it is a weak solution satisfying 
\[a_j\in L^\infty([0,T); l^2)\cap L^2([0,T); H^1) \ \ \ \forall \ \ j\geq 0\]
and 
\begin{equation}\notag
\|a(t)\|^2_{l^2}+2\nu\int_0^t \|a(\tau)\|^2_{H^1}\, d\tau\leq \|a(0)\|^2_{l^2}+\sum_{j=0}^\infty \int_0^tf_j(\tau) a_j(\tau)\, d\tau
\end{equation}
for all $t\in[0,T]$.
\end{Definition}
In \cite{Fi}, the author showed uniqueness of Leray-Hopf solution for $d\geq 1$, and uniqueness in the case of $d< 1$ but with additional decay condition. Namely,

\begin{Theorem}\label{thm-unique1}
Let $f_j\equiv 0$ for all $j\geq 0$ and $\{a_j(0)\}\in l^2$. If either 
\[d\geq 1\] or 
\[d <1, \ \ \ \mbox{and} \ \ \ a_j=o(\lambda_j^{\frac{d-1}{2}}) \ \ \mbox{as} \ \ j\to\infty,\]
then there exists a unique Leray-Hopf solution to (\ref{nse-force}) with initial data $\{a_j(0)\}$.
\end{Theorem}
We point out that in the case of $d<1$ the decay condition $a_j=o(\lambda_j^{\frac{d-1}{2}})$ implies that the solution is in $H^s$ for any $s<\frac{1-d}{2}$. 

With appropriate non-zero external forcing, Filonov and Khodunov \cite{FK} proved that uniqueness also holds for (\ref{nse-force}) with $d\geq1$. 

\begin{Theorem}\label{thm-unique2}
Let $d\geq1$ and $\{a_j(0)\}\in l^2$. If the forcing satisfies
\begin{equation}\label{force-h1}
\sum_{j=0}^\infty \lambda_j^{-2} \int_0^T f_j^2(t)\, dt<\infty,
\end{equation}
there exists a unique Leray-Hopf solution to (\ref{nse-force}) with initial data $\{a_j(0)\}$.
\end{Theorem}

Moreover they showed the opposite for $d<1$. 
\begin{Theorem}\label{thm-non-unique1}
Let $d<1$ and $a_j(0)=0$ for all $j\geq 0$. There exist a time $T>0$ and a sequence of functions $\{f_j(t)\}$ satisfying (\ref{force-h1}) such that the system (\ref{nse-force}) with initial data $\{a_j(0)\}$ has two different Leray-Hopf solutions.
\end{Theorem}

Thus $d=1$ is the sharp threshold to distinguish the uniqueness regime from the non-uniqueness regime. The presence of the external forcing is crucial for the non-uniqueness construction.  The idea of the construction is explained briefly as follows. Define the linear operator 
\[A(\{a_j\})=\{\lambda_j^2a_j\}\]
and bilinear operator 
\[B(\{a_j\}, \{v_j\})=\left\{-\lambda_j^{\theta} a_{j-1}v_{j-1}+\lambda_{j+1}^\theta a_j v_{j+1}\right\}.\]
Suppose system (\ref{nse-force}) has two solutions $a^+$ and $a^-$. Denote $v=\frac12(a^++a^-)$ and $g=\frac12(a^+-a^-)$. Then $(v,g)$ satisfies
\begin{equation}\label{eq-vg}
\begin{split}
\frac{d}{dt}v+Av+B(v,v)+B(g,g)=&\ f,\\
\frac{d}{dt}g+Ag+B(v,g)+B(g,v)=&\ 0
\end{split}
\end{equation}
with $v(0)=a(0)$ and $g(0)=0$. The authors constructed $v(t)$ and non-vanishing $g(t)$ which satisfy the $g$ equation of (\ref{eq-vg}). Substituting $v(t)$ and $g(t)$ into the $v$ equation of (\ref{eq-vg}) specifies the external forcing $f(t)$. The construction guarantees $f$ satisfies (\ref{force-h1}) and hence the existence of a Leray-Hopf solution is ensured. 

Interestingly, there is a similarity between this non-uniqueness result for the dyadic NSE model and a recent result for the actual 3D NSE with external forcing in \cite{ABC} by Albritton, Bru\'e and Colombo.  The authors constructed two Leray-Hopf solutions for the forced NSE in $\mathbb R^3$ with zero initial data

\begin{equation}\label{nse-3d-force}
\begin{split}
u_t+(u\cdot\nabla) u+\nabla p=&\  \Delta u+f,\\
\nabla\cdot u=&\ 0.
\end{split}
\end{equation}
Namely, they showed

\begin{Theorem}\label{thm-ABC}
There exists $T>0$ and $f\in L^1(0,T; L^2(\mathbb R^3_+))$ such that there are two distinct Leray-Hopf solutions to the NSE (\ref{nse-3d-force}) with the same forcing $f$ and initial data $u_0\equiv 0$.
\end{Theorem}
The result was cleverly built on the instability and non-uniqeness works \cite{Vishik1, Vishik2} for the Euler equation by Vishik. Among the two solutions, one is an unstable background solution similar to the unstable 2D vortex constructed in \cite{Vishik1, Vishik2}, and another is a trajectory on the unstable manifold associated to the background solution. In particular, the forcing $f$ has the form 
\[f(x,t)=t^{-\frac32} F\]
for a smooth and compactly supported profile function $F$. One notices that the forcing $f$ blows up as time approaches zero. Coincidently, the forcing in Theorem \ref{thm-non-unique1} for the dyadic NSE has the property
\[\lim_{t\to0}\lim_{j\to\infty} f_j(t) =\infty\]
which agrees with the feature of the forcing in the non-uniqueness result for the 3D NSE in Theorem \ref{thm-ABC}.  


\medskip

\subsection{Vanishing viscosity limit and dissipation anomaly}
\label{sec-attrator}

\noindent

Kolmogorov predicted energy cascade produces the so called dissipation anomaly, also referred to as Kolmogorov's zeroth law, which is
characterized by the persistence of non-vanishing energy dissipation in the limit of vanishing viscosity. This phenomenon is related to anomalous dissipation, namely the failure of the energy to be conserved despite the absence of viscosity discussed in Subsection~\ref{sec-open}. 

Rigorous justification of Kolmogorov's zeroth law and phenomenological scaling laws (which remain open up to this day for the NSE itself) were also realized for the dyadic NSE model. In the context of this model the properties of fixed points played an important role in establishing the results. 
We briefly review these works below. 

Cheskidov and Friedlander \cite{CF} showed dissipation anomaly through the study of the global attractor, which is a unique fixed point. Moreover the authors proved that the global attractor of the dyadic NSE (\ref{nse-force}) converges to the global attractor of the dyadic Euler model (\ref{euler-force}) in the vanishing viscosity limit. Again, the persistence of positivity is a crucial property in the analysis. 

\begin{Theorem}\label{thm-CF}\cite{CF}
Let $\nu>0$ in (\ref{nse-force}). Then
\begin{itemize}
\item [(i)] For each $\nu>0$, there is a unique fixed point $a^{\nu, *}$. The fixed point is a global attractor. 
\item [(ii)] Let $a^*$ be the fixed point in Theorem \ref{thm-CFP}. We have
\[a^{\nu, *} \to a^{*}, \ \ \mbox{as} \ \ \nu\to0.\]
\item [(iii)] Denote the energy dissipation rate for (\ref{nse-force}) by 
\[\varepsilon^\nu=\nu \|a^{\nu, *}\|^2_{H^1}=\left<f, a^{\nu, *}\right>.\]
Let $\varepsilon=a^*_0 f_0$ be the energy dissipation rate for the inviscid case as in Theorem \ref{thm-CFP}. We have
\[\lim_{\nu\to0} \varepsilon^{\nu}=\varepsilon>0.\]
\end{itemize}
\end{Theorem}
Thus Kolmogorov's zeroth law is proved for the dyadic NSE model. \\



\medskip

\subsection{Tao's averaged 3D NSE based on dyadic models} \label{subs:AveragedNSE}

\noindent

To understand the supercriticality barrier for the global regularity problem for the 3D NSE, Tao \cite{Tao} suggested a modified NSE with an averaged nonlinear structure and showed finite time blow-up for this equation. 

This modified 3D NSE is built on an ingeniously designed dyadic-type model where the energy does not spread as it cascades to high modes. As discussed in Subsection~\ref{subs:RegBlow}, solutions to the classical viscous dyadic model blow up when $d<-1$ \cite{Ch}, which is outside of physically relevant range in 3D. However blow-up occurs in the five-dimensional case.
As Tao remarked in his blog \cite{Tao-blog}, this is a common threshold for models where the blow-up is known. In fact, solutions of the viscous dyadic model are regular in the whole range of the intermittency dimension  $d\in[0,1]$ due to the result by Barbato, Morandin and Romito \cite{BMR}. The regularization effect comes from the nonlinear term as it spreads the energy as it cascades to high modes and prevents a blow-up.

To achieve a blow-up Tao has created a dyadic-type model with an added delay mechanism that forces an almost complete transfer of the energy from $a_j$ to $a_{j+1} $ before it starts moving to $a_{j+2}$, and hence destroys the regularizing effect of the nonlinear term described in Subsection~\ref{sebs:regularizing}. Based on this model, Tao proposed the approximating NSE
\begin{equation}\label{nse-tao}
u_t=\Delta u+\widetilde {\mathcal B}(u,u)
\end{equation}
where the bilinear operator $\widetilde {\mathcal B}$ is an averaged version of the bilinear operator 
\[\mathcal B(u,v)=-\frac12\mathbb P\left((u\cdot\nabla)v+(v\cdot\nabla) u\right)\]
taking the form
\begin{equation}\notag
\widetilde {\mathcal B}(u,v)=\int_{\Omega} m_{3,\omega}(D)Rot_{3,\omega}B\left(m_{1,\omega}(D)Rot_{1,\omega}, m_{2,\omega}(D)Rot_{2,\omega}\right)\, d\mu(\omega)
\end{equation}
for some probability space $(\Omega, \mu)$, some spatial rotation operators $Rot_{i,\omega}$ for $i=1,2,3$, and some Fourier multipliers $m_{i,\omega}$ of order 0. The operator $\widetilde {\mathcal B}$ obeys almost all of the estimates satisfied by the original operator $\mathcal B$.
In addition, the cancellation law
\[\int \widetilde {\mathcal B}(u,u)\cdot u\, dx=0 \]
holds. 
It was shown in \cite{Tao} that

\begin{Theorem}\label{thm-tao}
There exists a Schwartz vector field $u_0$ with $\nabla\cdot u_0=0$ such that there is no global-in-time mild solution to the averaged NSE (\ref{nse-tao}) with initial data $u_0$.
\end{Theorem}


The construction of the averaged bilinear operator $\widetilde {\mathcal B}$ only retains some nonlinear interactions from the original operator $\mathcal B$ and is designed with appropriate weights on the interactions to generate a specific blow-up mechanism. The nonlinearities of $\widetilde {\mathcal B}$ are a finite linear combination of local cascade operators implemented with a delay mechanism. In particular, the author called it a mechanism of a von Neumann machine (c.f. \cite{Tao-blog}). This result indicates that attempts to show global regularity for the 3D NSE using the abstract structure of the operator $\mathcal B$ is not likely to succeed. Instead, one needs to investigate the finer structure of the nonlinear term and use a dynamical approach. 

\medskip

\subsection{Ladyzhenskaya's modified NSE and the corresponding dyadic model} 

\noindent

The classical NSE describes the motion of an incompressible Newtonian flow in which the dissipative effect is captured by a linear relation between the Cauchy stress tensor and rate of strain tensor. To describe the isotropic turbulent flows with Kolmogorov's similarity hypothesis \cite{K41}, Smagorinsky \cite{Sma1, Sma2} proposed a class of modified NSEs in the meteorological context where nonlinear lateral diffusions were formulated.  Roughly speaking, the nonlinear diffusions have stronger regularization effect compared to that of the standard NSE. Later on, a modified NSE with more general nonlinear dissipation was suggested by Ladyzhenskaya \cite{Lady7, Lady8, Lady4}. 
The Ladyzhenskaya modified NSE reads as
\begin{equation}\label{nse-lady}
\begin{split}
u_t+(u\cdot\nabla) u+\nabla p=&\ \div (T(D))+f,\\
\nabla\cdot u=&\ 0
\end{split}
\end{equation}
where $D$ is the symmetric part of the velocity gradient tensor given by
\begin{equation}\notag
D=\frac12\left(\nabla u+(\nabla u)^T\right)
\end{equation}
and the stress tensor $T$ is a function of $D$. Specifically, $T$ satisfies the conditions that the entries $T_{ik}$ are continuous functions of $\frac{\partial u_i}{\partial x_k}$ and for some $\mu>0$
\begin{equation}\notag
\begin{split}
|T_{ik}(D)|\leq&\ c_1\left(1+|D|^{2\mu}\right) |D|, \\
T_{ik}(D)\left(\frac{\partial u_i}{\partial x_k}\right)\geq &\ \nu_0 D^2+\nu_1 D^{2+2\mu},\\
\int_{\Omega} \left(T_{ik}(D')-T_{ik}(D'')\right)\left(\frac{\partial u_i'}{\partial x_k} -\frac{\partial u_i''}{\partial x_k}\right)\, dx
\geq&\ \nu_2\int_{\Omega} \sum_{i,k=1}^3\left(\frac{\partial u_i'}{\partial x_k} -\frac{\partial u_i''}{\partial x_k}\right)^2\, dx
\end{split}
\end{equation}
where $u'$ and $u''$ are arbitrary smooth divergence free vectors and equal on $\partial \Omega$, and $c_1, \nu_0, \nu_1, \nu_2$ are constants. It was pointed out in \cite{Lady5} that the conditions on $T$ are physical for models of particle collision, such as in the study of Boltzmann equation. For $\mu\geq \frac14$, Ladyzhenskaya \cite{Lady4} showed existence of a unique global solution to (\ref{nse-lady}) in the energy space $L^2$. In contrast, the global solvability for the 3D standard NSE is still open. Other important properties of this modified NSE, such as the existence of a compact global attractor, were established by Ladyzhenskaya and her collaborators. Expository background on the modified NSEs can be found in \cite{Nec, Sma3}.  

In order to address the questions of well-posedness and partial regularity for general value $\mu\geq 0$, Friedlander and Pavlovi\'c \cite{FP2} proposed a dyadic model for the modified NSE (\ref{nse-lady}) as follows
\begin{equation}\label{FP-s2}
\frac{d}{dt} a_j=\lambda_j^{\frac52}a_{j-1}^2-\lambda_{j+1}^{\frac52}a_ja_{j+1}-\nu \lambda_j^{5\mu+2} a_j^{2\mu+1}, \ \ \ j\geq 0
\end{equation}
with $a_{-1}=0$. Existence of weak solutions of (\ref{FP-s2}) is proved for $\mu>\frac{1}{10}$. Namely,
\begin{Theorem}\label{thm-FP-modified}
Let $\mu>\frac1{10}$. If the initial data $a_j(0)$ satisfies certain decay condition, there exist weak solutions to (\ref{FP-s2}) with the initial data $a_j(0)$.
\end{Theorem}
 
Moreover the authors showed an upper bound for the Hausdorff dimension of the potential singular set for (\ref{FP-s2}) with $\mu<\frac12$.

\bigskip

\section{Results of dyadic NSE and Euler models with stochastic forcing}
Kolmogorov's phenomenological theory of turbulence relies heavily on statistical analysis and approaches. There is an extensive literature on the study of NSE and Euler equation with stochastic forcing. In the context of dyadic models, stochastically forced models have also been considered by many authors. Some of Kolmogorov's laws are rigorously justified for these dyadic models with certain stochastic forcing. 

\subsection{Additive noise on low modes}

\noindent

We consider the following model with additive noise,
\begin{equation}\label{nse-ss}
\begin{split}
d a_0+(\nu a_0+a_0a_1)\, dt=&\ \sigma \,dW,\\
\frac{d}{dt}a_j+\nu \lambda_j^2a_j+\lambda_{j+1}^\theta a_ja_{j+1}-\lambda_j^{\theta}a^2_{j-1}=&\ 0, \ \ j\geq 1
\end{split}
\end{equation}
where $W$ is a Brownian motion and constant $\sigma$ denotes the intensity of the noise. Note that the stochastic forcing only acts on a large scale. Denote 
\[S=(\Omega, \mathcal F, \{\mathcal F_t\}, \mathbb P, W)\] by the stochastic basis. 

We first recall the existence and uniqueness results for (\ref{nse-ss})
established by Friedlander, Glatt-Holtz and Vicol \cite{FGV}.
\begin{Theorem}\label{thm-ss-ex}
Let $\nu>0$ and $a^0\in l^2$. 
\begin{itemize}
\item [(i)] For $\theta\geq 1$ there exists a martingale solution $(a, S)$ to (\ref{nse-ss}) with initial data $a^0$ which satisfies
\begin{equation}\notag
\begin{split}
a\in&\ L^2\left( \Omega; L^\infty ([0,T]; l^2)\cap L^2([0,T]; H^1) \right) \ \ \ \mbox{for any} \ \ T>0,\\
a_j\in&\ C([0,\infty)) \ \ \mbox{a.s. for each} \ \ j\geq 0.
\end{split}
\end{equation}
\item [(ii)] If the components of the initial data $a^0$ are positive, then for any martingale solution $(a, S)$ we have $a_j(t)>0$ for all $j\geq 0$ and $t\geq 0$. Moreover, there is a realization $(u, S)$ such that 
\begin{equation}\label{moment}
\mathbb E |a(t)|^2+2\nu \int_0^t\mathbb E|a(s)|^2_{H^1}\, ds \leq |a^0|^2 +t\sigma^2.
\end{equation}
\item [(iii)] Let $\nu>0$ and $\theta\geq 1$. There exists a stationary martingale solution $(\bar a^{\nu}, S)$ satisfying (\ref{moment}).
\item [(iv)] Let $\theta\in[1,2]$ and fix a stochastic basis $S$. There exists a unique pathwise solution $a=a(t, a^0, W)$ of (\ref{nse-ss}) satisfying (\ref{moment}) with an equality. Moreover, $a(t, a^0, W)$ depends continuously on $a^0\in l^2$ and $W\in C([0,T])$. 
 \end{itemize}
\end{Theorem}
Recall $\theta=\frac{5-d}{2}$. Hence the uniqueness result in Theorem \ref{thm-ss-ex} (iv) is consistent with the uniqueness result in the deterministic case, see Theorem \ref{thm-unique2}.

Moreover the authors gave a rigorous proof for some of Kolmogorov's laws including the zeroth law - anomalous dissipation and the dissipation anomaly, through analysis of the statistically stationary solutions. 

\begin{Theorem}\label{thm-ss-steady} Let $\nu=0$.
\begin{itemize}
\item [(i)]
There exists a stationary martingale solution $(\bar a, S)$ to (\ref{nse-ss}), satisfying 
\begin{equation}\notag
\begin{split}
\bar a\in&\ L^\infty_{loc}([0,\infty); H^s) \ \ \ \mbox{for any} \ \ s<\frac13\theta,\\
\bar a_j\in&\ C([0,\infty)) \ \ \ \mbox{for each} \ \ j\geq 0, a.s.
\end{split}
\end{equation}
and 
\begin{equation}\label{ss-spectrum}
\sup_{j\geq 0}\lambda_j^{\frac23\theta} \mathbb E(\bar a_j^2)\leq C\sigma^{\frac43}
\end{equation}
for an absolute constant $C>0$. 
\item [(ii)]
Such solution $\bar a$ from (i) may be obtained as an inviscid limit, i.e. there exists Borel probability measures $\{\mu_{\nu_k}\}$ and $\mu_0$ on $l^2$ such that
\begin{equation}\notag
\mu_{\nu_k} \to \mu_0 \ \ \mbox{weakly in} \ \ H^{-\frac12} \ \ \mbox{as} \ \ \nu_k\to 0,
\end{equation}
with $\mu_{\nu_k}(\cdot)= \mathbb P(\bar a^{\nu_k}\in \cdot)$.
\item [(iii)]
The inviscid stationary solutions $\bar a$ has a constant mean energy flux
\begin{equation}\label{anol}
\varepsilon:=\mathbb E(\lambda_j^\theta \bar a_j^2a_{j+1}) =\frac12\sigma^2>0, \ \ \forall j\geq 0,
\end{equation}
and 
\begin{equation}\label{onsager}
\lim_{j\to\infty} \lambda_j^\theta \mathbb E |\bar a_j|^3>0.
\end{equation}
\end{itemize}
\end{Theorem} 
One can see that (\ref{anol}) indicates that the stationary inviscid solutions $\bar a$ realizes the anomalous dissipation of energy. Moreover, combining (\ref{ss-spectrum}) and (\ref{anol}) we infer the energy spectrum has the scaling
$\varepsilon^{\frac23}\lambda_j^{-\frac23\theta-1}$ which is consistent with Kolmogorov's spectrum scaling law when $\theta=1$ and hence the intermittency dimension $d=3$, which is referred as Kolmogorov's regime. In the end, it follows from (\ref{onsager}) that the inviscid steady state $\bar a$ has regularity below Onsager's critical space. 
 
The uniqueness and attraction properties of the invariant measure were also investigated in \cite{FGV}. Beside their own importance, these properties play a crucial role to show dissipation anomaly in the inviscid limit, in the agreement with the situation of the deterministic setting (c.f. \cite{CF}). We state the results below.

\begin{Theorem}\label{thm-invariant}
Consider (\ref{nse-ss}) with $\nu>0$, $\theta\in[1,2)$ and initial data $a^0\in l^2$. There exists a unique invariant measure $\mu_{\nu}$ of the corresponding Markov semigroup which is ergodic. Moreover, the invariant measure $\mu_{\nu}$ obeys the attraction properties of mixing, strong law of large numbers and central limit. 
\end{Theorem}

\begin{Theorem}\label{thm-ss-anomaly}
Let $a^{\nu}(\cdot, a^0)$ be the unique solution of (\ref{nse-ss}) with $\theta\in[1,2)$ and initial data $a^0\in l^2$. We have
\begin{equation}\notag
\begin{split}
&\lim_{\nu\to 0}\lim_{T\to\infty}\nu \mathbb E|a^{\nu}(T, a^0)|^2_{H^1}=\frac12\sigma^2,\\
&\lim_{\nu\to 0}\lim_{T\to\infty}\frac{\nu}{T} \int_0^T|a^{\nu}(t, a^0)|^2_{H^1}\, dt=\frac12\sigma^2, a.s. \ \ \mbox{for any} \ \ a^0\in l^2. 
\end{split}
\end{equation}
\end{Theorem}
The anomalous dissipation rate $\frac12\sigma^2$ in Theorem \ref{thm-ss-steady} is the same as the rate of dissipation anomaly in Theorem \ref{thm-ss-anomaly}.


\subsection{Additive noise on all modes}

\noindent

The model with additive noise on all scales 
\begin{equation}\label{nse-ss2}
da_j +\left(\nu\lambda_j^2a_j+\lambda_{j+1}^\theta a_ja_{j+1}-\lambda_j^\theta a_{j-1}^2\right)\, dt=\sigma_j dW_j
\end{equation}
was studied by Romito \cite{Rom}. Here, $W_j$ is one-dimensional Brownian motion and $\sigma_j\in \mathbb R$ for any $j\geq 0$. The author showed path-wise uniqueness for (\ref{nse-ss2}) with $\theta\in(2, 3]$ and finite time blow-up with positive probability when $\theta>3$ with certain assumption on the intensity of the noise.  

\begin{Theorem}\label{thm-Rom1}
Let $\theta\in(2,3]$. Assume there exists a constant $\alpha_0>\max\{\frac12(\theta-3), \theta-3\}$ such that
\begin{equation}\label{assu-rom}
\sup_{j\geq 0} \lambda_j^{\alpha_0} \sigma_j<\infty.
\end{equation}
There exists a path-wise unique solution of (\ref{nse-ss2}) with initial data in $H^{\theta-2}$ in the class of Galerkin martingale solutions. 
\end{Theorem}

\begin{Theorem}\label{thm-Rom2}
Let $\theta>3$ and (\ref{assu-rom}) holds. There exists initial data $a^0$ such that for every energy martingale weak solution starting at $a^0$, the probability of having a finite stopping time is positive. 
\end{Theorem}

\subsection{Multiplicative noise}

\noindent

Stochastic dyadic models with multiplicative noise have also been investigated. In particular, Barbato, Flandoli and Morandin \cite{BFM2} considered the inviscid model with multiplicative stochastic forcing in Stratonovich form
\begin{equation}\label{nse-ss3}
da_j+\left(\lambda_{j+1} a_ja_{j+1} -\lambda_j a_{j-1}^2\right)\, dt=\sigma\lambda_j a_{j-1}\circ dW_{j-1}-\sigma\lambda_{j+1} a_{j+1}\circ dW_j
\end{equation}
with Brownian motions $\{W_j\}$. Denote the total energy by $\mathcal E (t)=\frac12\sum_{j=0}^\infty a_j^2(t)$. The authors showed anomalous dissipation in the following sense. 

\begin{Theorem}\label{thm-BFM1}
Let $a(t)$ be the unique solution of (\ref{nse-ss3}) with initial data $a^0\in l^2$ which satisfies $\mathcal E (t)<\infty$ for all $t>0$. Then we have
\begin{equation}\notag
\mathbb P(\mathcal E(t)=\mathcal E(0))<1, \ \ \forall \ \ t>0.
\end{equation}
Moreover, for any $\varepsilon>0$ there exists $t$ such that
\begin{equation}\notag
\mathbb P(\mathcal E(t)<\varepsilon) >0.
\end{equation}
\end{Theorem}
Despite the formal conservation of the energy for system (\ref{nse-ss3}), such anomalous dissipation arises from the forward energy cascade mechanism, just like the deterministic dyadic Euler model. Consequently, the authors further proved that global regular solutions can not exist by employing the stochastic tool of birth and death process. 
Later on, Barbato and Morandin \cite{BM2} extended the work of \cite{BFM2} to more general dyadic models forced by multiplicative noise, including the inviscid GOY and Sabra models. 

For the same inviscid model (\ref{nse-ss3}) with multiplicative noise, the authors showed weak uniqueness in the class of exponentially integrable solutions in \cite{BFM3}. Namely,

\begin{Theorem}\label{thm-BFM2}
Consider (\ref{nse-ss3}) with initial data $a^0\in l^2$. A weak solution $(\Omega, F_t, P, W, a)$ of (\ref{nse-ss3}) is said to be exponentially integrable on $[0,T]$ if
\begin{equation}\notag
E^P\left[e^{\frac1{\sigma^2}\int_0^T \sum_{j=0}^\infty a_j^2(t)\, dt} \left(1+\int_0^T a_i^4(t)\, dt\right)^2\right]<\infty, \ \ \forall \ \ i\geq 0.
\end{equation}
There is weak uniqueness in the class of exponentially integrable solutions for (\ref{nse-ss3}).
\end{Theorem}

Weak uniqueness is in the sense of uniqueness of the law of the process. Recall that in the deterministic case, there are non-unique weak solutions for the dyadic Euler model. Nevertheless, strong (path-wise) uniqueness for (\ref{nse-ss3}) remains an open question. As indicated by the authors, their approach of showing uniqueness builds on the Girsanov transformation which turns a nonlinear system to a linear system. This approach naturally leads to weak uniqueness. 



\bigskip

\section{Dyadic models for magnetohydrodynamics}

The PDEs for the incompressible magnetohydrodynamics (MHD) 
\begin{equation}\label{mhd}
\begin{split}
u_t+(u\cdot\nabla) u-(B\cdot\nabla) B+\nabla P=&\ \nu\Delta u,\\
B_t+(u\cdot\nabla)B-(B\cdot\nabla) u=&\ \mu\Delta B,\\
\nabla\cdot u=0, \ \ \nabla\cdot B=&\ 0
\end{split}
\end{equation}
describe electrically conducting fluids in geophysics and astrophysics. In system (\ref{mhd}), $u$, $B$ and $p$ are the unknowns, denoting respectively the velocity field, magnetic field and scalar pressure function. The constants $\nu$ and $\mu$ represent the kinetic viscosity and magnetic diffusivity respectively. Despite many similar features with the NSE, the MHD system (\ref{mhd}) is more challenging due to the intricate interactions of the velocity field and magnetic field. The open questions for the pure fluid equations discussed in Subsection \ref{sec-open} remain open for the MHD system. Inspired by the fruitful results for the dyadic models for the NSE and Euler equation, it is possible to advance the understanding of the MHD by studying reduced models. 
In particular, a class of dyadic models were introduced by Dai in \cite{Dai-20} including both of the Desnyansky-Novikov and Obukhov types of nonlinear structures. Particular case of these models are consistent with some dyadic models suggested by physicists, for instance, see the work \cite{PSF} of Plunian, Stepanov and Frick. We point out that the majority of the work from the physics community concerns numerical study of these dyadic models. 


We recall two particular models from \cite{Dai-20} which contain Desnyansky-Novikov nonlinear terms and exhibit different energy cascade mechanisms, 
\begin{equation}\label{mhd1}
\begin{split}
\frac{d}{dt} a_j+\nu\lambda_j^2 a_j+\lambda_{j}^\theta a_ja_{j+1}-\lambda_{j-1}^\theta a_{j-1}^2+\lambda_{j}^\theta b_jb_{j+1}-\lambda_{j-1}^\theta b_{j-1}^2=&\ f_j,\\
\frac{d}{dt} b_j+\mu\lambda_j^2 b_j-\lambda_{j}^\theta a_jb_{j+1}+\lambda_j^\theta b_ja_{j+1}=&\ 0
\end{split}
\end{equation}
for $j\geq 0$ with $a_{-1}=b_{-1}=0$ and external forcing $f=(f_0, f_1, f_2, ...)$, and
\begin{equation}\label{mhd2}
\begin{split}
\frac{d}{dt} a_j+\nu\lambda_j^2 a_j+\lambda_{j}^\theta a_ja_{j+1}-\lambda_{j-1}^\theta a_{j-1}^2-\lambda_{j}^\theta b_jb_{j+1}+\lambda_{j-1}^\theta b_{j-1}^2=&\ f_j,\\
\frac{d}{dt} b_j+\mu\lambda_j^2 b_j+\lambda_{j}^\theta a_jb_{j+1}-\lambda_j^\theta b_ja_{j+1}=&\ 0.
\end{split}
\end{equation}
Here $a_j(t)=\|u(t)\|_{L^2}$ and $b_j(t)=\|B(t)\|_{L^2}$, and $\theta=\frac{5-d}{2}$ (under the assumption that the intermittency dimension of the velocity field and magnetic field is the same). Denote the total energy by
\[E(t)=\frac12\sum_{j=0}^\infty \left(a_j^2(t)+b^2_j(t)\right)\]
and cross helicity 
\[H(t)=\sum_{j=0}^\infty a_j(t)b_j(t).\]
When $\nu=0$ and $f\equiv 0$, system (\ref{mhd1}) formally conserves the total energy, but not the cross helicity; while (\ref{mhd2}) conserves both the total energy and cross helicity. 
For positive $a_j$ and $b_j$, we note that there is only forward energy cascade in the velocity equations of model (\ref{mhd1}), while there are both forward and backward energy cascades in the velocity equations of (\ref{mhd2}). The energy transfer between velocity shells and magnetic field shells is rather complex. In particular positivity is not known to be preserved for the dyadic MHD models (\ref{mhd1}) and (\ref{mhd2}). This is one of the major differences from dyadic NSE models. Without the positivity property the actual energy transfer within the dyadic MHD models is substantially more complicated than that of the pure fluid models. 


We will briefly discuss certain results that were recently proved for the MHD models (\ref{mhd1}) and (\ref{mhd2}). Dai and Friedlander \cite{DF1}
considered (\ref{mhd1}) with $\nu=\mu=0$ and forcing acting only on large scales. Without loss of generality, we assume $f_0\geq 0$ and $f_i\equiv 0$ for $j\geq 1$. 

\begin{Theorem}\label{thm-DF1}
Let $\theta>0$. The solution $(a(t), b(t))$ of (\ref{mhd1}) with positive initial data develops blow-up at a finite time in the $H^s$ norm with $s>\frac13\theta$.
\end{Theorem}
The proof does not require the solution to stay positive, although it is necessary to start with positive initial data to produce the finite time blow-up. In contrast, the proof breaks down to show finite time blow-up for system (\ref{mhd2}) because the bi-directional energy cascade mechanism is more intricate. 

In \cite{DF1} the authors also studied the steady states of (\ref{mhd1}) and (\ref{mhd2}) with $\nu=\mu=0$ and forcing $f_0>0$. The steady states of (\ref{mhd1}) satisfy the exact form
\begin{equation}\label{fixed1}
a^*_j=A_0\lambda^{\frac16\theta} f_0^{\frac12} \lambda_j^{-\frac13\theta}, \ \ b^*_j=B_0\lambda^{\frac16\theta} f_0^{\frac12} \lambda_j^{-\frac13\theta}
\end{equation}
with $A_0^2+B_0^2=1$. The steady states of (\ref{mhd2}) also have the form (\ref{fixed1}) with the free parameters $A_0$ and $B_0$ satisfying $A_0^2-B_0^2=1$. Obviously, in contrast with the dyadic Euler model which has a unique fixed point, the steady states of (\ref{mhd1}) and (\ref{mhd2}) form a set of infinitely many points with two degrees of freedom. In general, it is challenging to study any possible attraction properties of such fixed point sets. Indeed, in the special case of $A_0=1$ and $B_0=0$, the fixed point is shown to be linearly unstable.

 The viscous case of $\nu>0$ and $\mu>0$ was studied by Dai and Friedlander in \cite{DF2}. The authors addressed the question of uniqueness and non-uniqueness of Leray-Hopf solutions and proved that

\begin{Theorem}\label{thm-DF2}
Let $\theta>0$, $a^0\in l^2$ and $b^0\in l^2$. Assume 
\begin{equation}\label{force-DF}
\sum_{j=0}^\infty \lambda_j^{-2} \int_0^T f_j^2(t)\, dt<\infty, \ \ \ \mbox{for any} \ \ T>0.
\end{equation}
There exists a Leray-Hopf solution to (\ref{mhd1}) (and (\ref{mhd2})) with the initial data $(a^0, b^0)$. 
\end{Theorem}

\begin{Theorem}\label{thm-DF3}
Let $0<\theta\leq 2$, $a^0\in l^2$, $b^0\in l^2$ and the forcing $f$ satisfy (\ref{force-DF}).
The Leray-Hopf solution to (\ref{mhd1}) (and (\ref{mhd2})) with the initial data $(a^0, b^0)$ is unique. 
\end{Theorem}

\begin{Theorem}\label{thm-DF4}
Let $\theta> 2$, $a^0=b^0=0$. There exist a time $T>0$ and functions $\{f_j(t)\}$ satisfying (\ref{force-DF}) such that 
the system (\ref{mhd1}) (and (\ref{mhd2})) with the initial data $(a^0, b^0)$ has at least two Leray-Hopf solutions $(a(t), b(t))$ with non-vanishing $b(t)$ on $[0,T]$.
\end{Theorem}
We note that the persistence of positivity is not necessary to establish the results of Theorems \ref{thm-DF2}, \ref{thm-DF3} and \ref{thm-DF4}. The non-uniqueness construction follows the approach of Filonov and Khodunov \cite{FK} for the pure fluid dyadic model which is discussed in Theorem \ref{thm-non-unique1} of Subsection \ref{subsec-unique}.
 


\medskip

\bigskip

\section*{Acknowledgement}

A. Cheskidov is partially supported by the NSF grant DMS-1909849.
M. Dai is partially supported by the NSF grants DMS-1815069 and DMS-2009422, and the AMS Centennial Fellowship. S. Friedlander is partially supported by the NSF grant DMS-1613135. A. Cheskidov and M. Dai are grateful to IAS for its hospitality in 2021-2022. They are also grateful for the hospitality of Princeton University in 2022-2023.

\bigskip



\begin{thebibliography}{XX}




\bibitem{ABC}
D. Albritton, E. Bru\'e and M. Colombo.
\newblock {\em Non-uniqueness of Leray solutions of the forced Navier-Stokes equations}.
\newblock Annals of Mathematics, Vol. 196: 415--455, 2022.

\bibitem{BFM}
D. Barbato, F. Flandoli, and F. Morandin.
\newblock {\em Energy dissipation and self-similar solutions for an unforced inviscid dyadic model}.
\newblock Trans. Amer. Math. Soc., 363 (4): 1925--1946, 2011.

\bibitem{BFM2}
D. Barbato, F. Flandoli, and F. Morandin.
\newblock {\em Anomalous dissipation in a stochastic inviscid dyadic model}.
\newblock Annals of Applied Probability, 21(6): 2424--2446, 2011.

\bibitem{BFM3}
D. Barbato, F. Flandoli, and F. Morandin.
\newblock {\em Uniqueness for a stochastic inviscid dyadic model}.
\newblock Proceedings of the American Mathematical Society, 138(7): 2607--2617, 2010.

\bibitem{BM}
D. Barbato and F. Morandin.
\newblock {\em Positive and non-positive solutions for an inviscid dyadic model: well-posedness and regularity}.
\newblock Nonlinear Differential Equations Appl., 20 (3): 1105--1123, 2013.

\bibitem{BM2}
D. Barbato and F. Morandin.
\newblock {\em Stochastic inviscid shell models: well-posedness and anomalous dissipation}.
\newblock Nonlinearity, 26 (7): 1919--1943, 2013.

\bibitem{BMR3}
D. Barbato, F. Morandin, and M. Romito.
\newblock {\em Global regularity for a logarithmically supercritical hyperdissipative dyadic equation}.
\newblock Dynamics of PDE, Vol. 11(1): 39-52, 2014.

\bibitem{BMR2}
D. Barbato, F. Morandin, and M. Romito.
\newblock {\em Global regularity for a slightly supercritical hyperdissipative Navier-Stokes system}.
\newblock Analysis and PDE, Vol. 7, No.8, 2014.

\bibitem{BMR}
D. Barbato, F. Morandin, and M. Romito.
\newblock {\em Smooth solutions for the dyadic model}.
\newblock Nonlinearity, 24 (11): 3083--3097, 2011.

\bibitem{BBV}
R. Beekie, T. Buckmaster, and V. Vicol.
\newblock {\em Weak solutions of ideal MHD which do not conserve magnetic helicity}.
\newblock Ann. of PDE, 6 (1), https://doi.org/10.1007/s40818-020-0076-1, 2020.

\bibitem{BF}
H. Bessaih and B. Ferrario.
\newblock {\em Invariant Gibbs measures of the energy for shell models of turbulence: the inviscid and viscous cases}.
\newblock Nonlinearity, 25(4): 1075--1097, 2012.


\bibitem{Bian}
L.A. Bianchi.
\newblock {\em Uniqueness for an inviscid stochastic dyadic model on a tree}.
\newblock Electronic Communications in Probability, 18: 1--12, 2013.

\bibitem{Bif}
L. Biferale.
\newblock {\em Shell models of energy cascade in turbulence}.
\newblock Annu. Rev. Fluid Mech., 35: 441468, 2003.




\bibitem{BJPV}
T. Bohr, M.H. Jensen, G. Paladin, and A. Vulpiani.
\newblock {\em Dynamical Systems Approach to Turbulence}.
\newblock Cambridge University Press, 1998.






\bibitem{BV}
T. Buckmaster, and V. Vicol.
\newblock {\em Nonuniqueness of weak solutions to the Navier-Stokes equation}.
\newblock Ann. of Math., 189(1):101--144, 2019.









\bibitem{Ch}
A. Cheskidov.
\newblock {\em Blow-up in finite time for the dyadic model of the Navier-Stokes equations }.
\newblock Trans. Amer. Math. Soc., 360 (10): 5101-5120, 2008.

\bibitem{CCFS}
A. Cheskidov, P. Constantin, S. Friedlander, and R. Shvydkoy.
\newblock {\em Energy conservation and Onsager's conjecture for the Euler equations}.
\newblock Nonlinearity, 21(6):1233--1252, 2008.

\bibitem{CD-discont}
A. Cheskidov and M. Dai.
\newblock {\em Discontinuity of weak solutions to the 3D NSE and MHD equations in critical and supercritical spaces}.
\newblock Journal of Mathematical Analysis and Applications, Vol. 481 (2), 123493, 2020.

\bibitem{CD-Kol}
A. Cheskidov and M. Dai.
\newblock {\em Kolmogorov's dissipation number and the number of degrees of freedom for the 3D Navier-Stokes equations}.
\newblock Proceedings of the Royal Society of Edinburg, Section A, Vol. 149, Issue 2: 429--446, 2019.

\bibitem{CD-norm}
A. Cheskidov and M. Dai.
\newblock {\em Norm inflation for generalized Navier-Stokes equations}.
\newblock Indiana Univ. Math. J., 63(3): 869--884, 2014.



\bibitem{CF}
A. Cheskidov and S. Friedlander.
\newblock {\em The vanishing viscosity limit for a dyadic model}.
\newblock Physica D, 238:783--787, 2009.

\bibitem{CFP1}
A. Cheskidov, S. Friedlander, and N. Pavlovi\'c.
\newblock {\em Inviscid dyadic model of turbulence: the fixed point and Onsager's conjecture}.
\newblock J. Math. Phys., 48 (6): 065503, 16, 2007.

\bibitem{CFP2}
A. Cheskidov, S. Friedlander, and N. Pavlovi\'c.
\newblock {\em An inviscid dyadic model of turbulence: the global attractor}.
\newblock Discrete Contin. Dyn. Syst., 26 (3): 781--794, 2010.

\bibitem{CL}
A. Cheskidov and X. Luo.
\newblock {\em Sharp nonuniqueness for the Navier-Stokes equations}.
\newblock Inventiones Mathematicae, Vol. 229: 987--1054, 2022.


\bibitem{CSreg}
A. Cheskidov and R. Shvydkoy.
\newblock {\em A unified approach to regularity problems for the 3D Navier-Stokes and Euler equations: the use of Kolmogorov's dissipation range}.
\newblock J. Math. Fluid Mech., Vol.16, Issue 2: 263--273, 2014.

\bibitem{CSint}
A. Cheskidov and R. Shvydkoy.
\newblock {\em Euler equations and turbulence: analytical approach to intermittency}.
\newblock SIAM J. Math. Anal., 46 (1): 353--374, 2014.



\bibitem{CZ}
A. Cheskidov and K. Zaya.
\newblock {\em Regularizing effect of the forward energy cascade in the inviscid dyadic model}.
\newblock Proc. Amer. Math. Soc., 144: 73--85, 2016.


\bibitem{CET}
P. Constantin, W. E, and E.Titi.
\newblock {\em Onsager's conjecture on the energy conservation for solutions of Euler's equation}.
\newblock Comm. Math. Phys., 165:207--209, 1994.

\bibitem{CLT}
P. Constantin, B. Levant, and E.Titi.
\newblock {\em Analytic study of the shell model of turbulence}.
\newblock Physica D: Nonlinear Phenomena, 219 (2): 120--141, 2006.



\bibitem{Dai-20}
M. Dai.
\newblock {\em Blow-up of a dyadic model with intermittency dependence for the Hall MHD}.
\newblock  Physica D: Nonlinear Phenomena, Vol. 428: 133066, 2021.







\bibitem{DF1}
M. Dai and S. Friedlander.
\newblock {\em Dyadic models for ideal MHD}.
\newblock Journal of Mathematical Fluid Mechanics, doi:10.1007/s00021-021-00640-9, 2021.

\bibitem{DF2}
M. Dai and S. Friedlander.
\newblock {\em Uniqueness and non-uniqueness results for dyadic MHD models}.
\newblock Journal of Nonlinear Science, https://doi.org/10.1007/s00332-022-09868-9, 2023.





\bibitem{DLS1}
C. De Lellis, and L. Sz\'ekelyhidi, Jr.
\newblock {\em Dissipative continuous Euler flows}.
\newblock Invent. Math., Vol.193 No. 2: 377--407, 2013.

\bibitem{DLS2}
C. De Lellis, and L. Sz\'ekelyhidi, Jr.
\newblock {\em Dissipative Euler flows and Onsager's conjecture}.
\newblock Journal of the European Mathematical Society,  16(7): 1467--1505, 2014.

\bibitem{DLS3}
C. De Lellis, and L. Sz\'ekelyhidi, Jr.
\newblock {\em The Euler equations as a differential inclusion}.
\newblock Ann. of Math.,  Vol.170 No.3: 1417--1436, 2009.


\bibitem{DN}
V. N.  Desnyansky and E. A. Novikov.
\newblock {\em Evolution of turbulence spectra toward a similarity regime}.
\newblock Izv. Akad. Nauk SSSR, Fiz. Atmos. Okeana, 10: 127--136, 1974.


\bibitem{DS}
E. I. Dinaburg and Y. G. Sinai.
\newblock {\em A quasi-linear approximation of three-dimensional Navier-Stokes system}.
\newblock Moscow Math. J., 1: 381--388, 2001.


\bibitem{ESS}
L. Escauriaza, G. Seregin, and \v Sver\'ak.
\newblock {\em $L^{3,\infty}$-solutions of Navier-Stokes equations and backward uniqueness}.
\newblock Uspekhi Mat. Nauk 58, 2(350):3--44, 2003.

\bibitem{Ey}
G. L. Eyink.
\newblock {\em Energy dissipation without viscosity in ideal hydrodynamics. I. Fourier analysis and local energy transfer}.
\newblock Phys. D, 78:222--240, 1994.

\bibitem{ES}
G. L. Eyink and K.R. Sreenivasan. 
\newblock {\em Onsager and the theory of hydrodynamic turbulence}.
\newblock Rev. Mod. Phys., 78, 2006.




\bibitem{Fi}
N. Filonov.
\newblock {\em Uniqueness of the Leray-Hopf solution for a dyadic model}.
\newblock Transactions of the American Mathematical Society, Vol. 369 (12): 8663--8684, 2017.

\bibitem{FK}
N. Filonov and P. Khodunov.
\newblock {\em Non-uniqueness of Leray-Hopf solutions for a dyadic model}.
\newblock St. Petersburg Math. J., Vol. 32: 371--387, 2021.

\bibitem{FGV}
S. Friedlander, N. Glatt-Holtz, and V. Vicol.
\newblock {\em Inviscid limits for a stochastically forced shell model of turbulent flow}.
\newblock Annales de l'Institut henri Poincar\'e - Probalilit\'es et Statistiques, 52(3): 1217--1247, 2016.

\bibitem{FP}
S. Friedlander and N. Pavlovi\'c.
\newblock {\em blow-up in a three-dimensional vector model for the Euler equations}.
\newblock Comm. Pure Appl. Math., 57 (6): 705--725, 2004.

\bibitem{FP2}
S. Friedlander and N. Pavlovi\'c.
\newblock {\em Remarks concerning modified Navier-Stokes equations}.
\newblock Discrete and Continuous Dynamical Systems, 10 (1-2): 269--288, 2004.

\bibitem{Fri}
U. Frisch.
\newblock {\em Turbulence: The Legacy of A. N. Kolmogrov}.
\newblock Cambridge University Press, Cambridge, 1995.


\bibitem{Gle}
E. B. Gledzer.
\newblock {\em System of hydrodynamic type admitting two quadratic integrals of motion}.
\newblock Soviet Phys. Dokl., 18: 216-217, 1973.

\bibitem{GLPG}
C. Gloaguen, J. L\'eorat, A. Pouquet and R. Grappin.
\newblock {\em A scalar model for MHD turbulence}.
\newblock Physica D.: Nonlinear Phenomena, 17(2):154--182, 1985.


\bibitem{Hopf}
E. Hopf.
\newblock {\em \"Uber die Anfangswertaufgabe f\"ur die hydrodynamischen Grundgleichungen}.
\newblock Math. Nachr., 4:213--231, 1951.


\bibitem{Is}
P. Isett.
\newblock {\em A Proof of Onsager's Conjecture}.
\newblock Ann. of Math., Vol.188 No.3: 1--93, 2018.

\bibitem{JL}
I. Jeong and D. Li.
\newblock {\em A blow-up result for dyadic models of the Euler equations}.
\newblock Communications in Mathematical Physics, 337:1027--1034, 2015.


\bibitem{JS1}
H. Jia and V. \v{S}ver\'ak.
\newblock {\em Are the incompressible 3d Navier-Stokes equations locally ill-posed in the natural energy space?}
\newblock J. Funct. Anal., Vol. 268(12): 3734--3766, 2015.


\bibitem{kaneda}
Y. Kaneda, T. Ishihara, M. Yokokawa, K. Itakura, and A. Uno. 
\newblock {\em Energy dissipation rate and energy spectrum in high resolution direct numerical simulations of turbulence in a periodic box}.
\newblock Physics of Fluids, 15 (2): 21--24, 2003.

\bibitem{KP-cheap}
N. Katz and N. Pavlovi\'c. 
\newblock {\em A cheap Caffarelli-Kohn-Nirenberg inequality for the Navier-Stokes equation with hyper-dissipation}.
\newblock Geom. Funct. Anal., 12 (2): 355--379, 2002.

\bibitem{KP}
N. Katz and N. Pavlovi\'c. 
\newblock {\em Finite time blow-up for a dyadic model of the Euler equations}.
\newblock Trans. Amer. Math. Soc., 357 (2): 695--708, 2005.

\bibitem{KZ}
A. Kiselev and A. Zlato\v{s}. 
\newblock {\em On discrete models of the Euler equation}.
\newblock Int. Math. Res. Not., 38: 2315--2339, 2005.

\bibitem{K41}
A. Kolmogorov.
\newblock The local structure of turbulence in incompressible viscous fluid for very large {R}eynold's numbers.
\newblock {\em C. R. (Doklady) Acad. Sci. URSS (N.S.)}, 30:301--305, 1941.

\bibitem{Lady-att1} 
O.A. Ladyzhenskaya.
\newblock {\em A dynamical system generated by the Navier-Stokes equations}.
\newblock Zap. Nau\v cn. Sem. Leningrad. Otdel. Mat. Inst. Steklov. (LOMI),  27: 91--115, 1972. English transl., J. Soviet. Math., 3:458--479, 1975.

\bibitem{Lady-att2} 
O.A. Ladyzhenskaya.
\newblock {\em Attractors for semigroups and evolution equatoins}.
\newblock Lezioni Lincei, 1988; Cambridge Univ. Press, Cambridge, 1991.

\bibitem{Lady0} 
O.A. Ladyzhenskaya.
\newblock {\em Solution ``in the large'' of boundary value problems for the Navier-Stokes equations in two space variables}.
\newblock Dokl. Akad. Nauk SSSR, 123:427--429, 1958. English transl., Soviet Phys. Dokl., 3:1128--1131, 1959; and Comm. Pure App. Math., 12:427--433, 1959.

\bibitem{Lady1} 
O.A. Ladyzhenskaya.
\newblock {\em Uniqueness and smoothness of generalized solutions of Navier-Stokes equations (Russian)}.
\newblock Zap. Nau\v cn. Sem. Leningrad. Otdel. Mat. Inst. Steklov. (LOMI),  5: 169--185, 1967.

\bibitem{Lady2} 
O.A. Ladyzhenskaya.
\newblock {\em New equations for the description of the motions of viscous incompressible fluids, and global solvability for their boundary value problems (Russian)}.
\newblock Trudy Mat. Inst. Steklov,  102: 85--104, 1967.

\bibitem{Lady3} 
O.A. Ladyzhenskaya.
\newblock {\em Certain nonlinear problems of the theory of continuous media (Russian)}.
\newblock In Proc. Internat. Congr. Math. (Moscow, 1966). Moscow: Izdat. ``Mir,''  : 560--573, 1968.

\bibitem{Lady7} 
O.A. Ladyzhenskaya.
\newblock {\em On nonlinear problems of continuum mechanics}.
\newblock Proc. Internat. Congr. Math. (Moscow, 1966),  Nauka, Moscow, 560--573, 1968.

\bibitem{Lady8} 
O.A. Ladyzhenskaya.
\newblock {\em On some new equations describing dynamics of incompressible fluids and on global solvability of boundary value problems to these equations}.
\newblock Trudy Steklov's Math. Institute,  102: 85--104, 1967.

\bibitem{Lady4} 
O.A. Ladyzhenskaya.
\newblock {\em On some modifications of the Navier-Stokes equations for large gradients of the velocities (Russian)}.
\newblock Zap. Nau\v cn. Sem. Leningrad. Otdel. Mat. Inst. Steklov. (LOMI),  7: 126--154, 1968.

\bibitem{Lady6} 
O.A. Ladyzhenskaya.
\newblock {\em Unique global solvability of the three-dimensional Cauchy problem for the Navier-Stokes equations in the presence of axial symmetry}.
\newblock Zap. Nau\v cn. Sem. Leningrad. Otdel. Mat. Inst. Steklov. (LOMI),  7: 155--177, 1968. English transl., Sem. Math. V.A. Steklov Math. Inst. Leningrad, 7:70--79, 1970.



\bibitem{Lady5} 
O.A. Ladyzhenskaya.
\newblock {\em The Mathematical Theory of Viscous Incompressible Flow}.
\newblock Second English ed., revised and enlarged. Translated from the Russian by Richard A. Silverman and John Chu. Mathematics and its Applications, Vol.2. New York-London-Paris: Gordon and Breach, Science Publishers, 1969.

\bibitem{LadyK} 
O.A. Ladyzhenskaya and A.A. Kiselev. 
\newblock {\em On the existence and uniqueness of the solution of the non-stationary problem for a viscous incompressible fluid}.
\newblock Izv. Akad. Nauk SSSR Ser. Mat. 21:665--680, 1957; English transl., Amer. Math. Soc. Transl. (2) 24:79--106, 1963.






\bibitem{Le}
J. Leray.
\newblock {\em Sur le mouvement d'un liquide visqueux emplissant l'espace}.
\newblock Acta Math., 63(1):193--248, 1934.


\bibitem{LP}
J.-L. Lions and G. Prodi.
\newblock {\em Un th\'eoreme d'existence et unicit\'e dans les \' equations de Navier-Stokes en dimension 2}.
\newblock C. R. Acad. Sci. Paris, 248:3519--3521, 1959.

\bibitem{Lorenz}
E.N. Lorenz.
\newblock {\em Deterministic nonperiodic flow}.
\newblock J. Atmos. Sci., 20:130--141, 1972.



\bibitem{LPPPV}
V.S. L'vov, E. Podivilov, A. Pomyalov, I. Procaccia, and D. Vandembroucq.
\newblock {\em Improved shell model of turbulence}.
\newblock Phys. Rev. E (3) 58: 1811--1822, 1998.

\bibitem{Mai0}
A.A. Mailybaev.
\newblock {\em Hidden scale invariance of intermittent turbulence in a shell model}.
\newblock Physical Review Fluids, 6, L012601, 2021.

\bibitem{Mai1}
A.A. Mailybaev.
\newblock {\em Shell model intermittency is the hidden self-similarity}.
\newblock Physical Review Fluids, 7, 034604, 2022.

\bibitem{Mai2}
A.A. Mailybaev.
\newblock {\em Solvable intermittent shell model of turbulence}.
\newblock Communications in Mathematical Physics, 388: 469--478, 2021.


\bibitem{MSV}
J.C. Mattingly, T. Suidan, and E. Vanden-Eijnden.
\newblock {\em Simple systems with anomalous dissipation and energy cascade}.
\newblock Communications in Mathematical Physics, 276(1): 189--220, 2007.




\bibitem{Nec}
J. Ne\v cas. 
\newblock {\em Theory of Multipolar Viscous Fluids}.
\newblock Academic Press, 1991.

\bibitem{Ob}
A. M. Obukhov. 
\newblock {\em Some general properties of equations describing the dynamics of the atmosphere}.
\newblock Izv. Akad. Nauk SSSR Ser. Fiz. Atmosfer. i Okeana, 7:695--704, 1971.



\bibitem{OY}
K. Ohkitani and M. Yamada.
\newblock {\em Temporal intermittency in the energy cascade process and local Lyapunov analysis in fully-developed model of turbulence}.
\newblock Progr. Theoret. Phys., 81: 329--341, 1989.

\bibitem{On}
L. Onsager.
\newblock {\em Statistical hydrodynamics}.
\newblock Nuovo Cimento (9), 6(Supplemento, 2(Convegno Internazionale di Meccanica Statistica)):279--287, 1949.

\bibitem{PSF}
F. Plunian, R. Stepanov and P. Frick.
\newblock {\em Shell models of magnetohydrodynamic turbulence}.
\newblock Physics Reports, vol. 523, 2013.


\bibitem{Prodi}
G. Prodi.
\newblock {\em Un teorema di unicita per el equazioni di Navier-Stokes}.
\newblock Ann. Mat. Pura Appl., 48: 173--182, 1959.

\bibitem{Rom}
M. Romito.
\newblock {\em Uniqueness and blow-up for a stochastic viscous dyadic model}.
\newblock Probability Theory and Related Fields, 158(3-4): 895--924, 2014.


\bibitem{Serrin}
J. Serrin.
\newblock {\em The initial value problem for the Navier-Stokes equations}.
\newblock Nonlinear Problems (R. Langer, ed.), pp. 69--98, The university of Wisconsin Press, Madison, 1963.




\bibitem{Sig}
E.D. Siggia.
\newblock {\em Model of intermittency in three-dimensional turbulence}.
\newblock Phys. Rev. A, 17: 1166--1176, 1978.

\bibitem{Sma1}
J. Smagorinsky.
\newblock {\em On the numerical integration of the primitive equations of motion for baroclinic flow in a closed region}.
\newblock Mon. Wea. Rev. 86, 3:457--466, 1958.

\bibitem{Sma2}
J. Smagorinsky.
\newblock {\em General circulation experiments with the primitive equations, Part I: The basic experiment}.
\newblock Mon. Wea. Rev. 91, 3:99--152, 1963.

\bibitem{Sma3}
J. Smagorinsky.
\newblock {\em Some historical remarks on the use of nonlinear viscosities, in ``Large eddy simulation of complex engineering and geophysical flows''}.
\newblock Cambridge University Press, 1993.

\bibitem{TT}
E. Tadmor and T. Tao.
\newblock {\em Velocity averaging, kinetic formulations, and regularizing effects in quasi-linear PDEs}.
\newblock Comm. Pure Appl. Math., Vol. 60(10): 1488--1521, 2007.

\bibitem{Tao}
T. Tao.
\newblock {\em Finite time blow-up for an averaged three-dimensional Navier-Stokes equation}.
\newblock J. Amer. Math. Soc., Vol. 29: 601--674, 2016.

\bibitem{Tao-blog}
T. Tao.
\newblock {\em Finite time blow-up for an averaged three-dimensional Navier-Stokes equation}.
\newblock T. Tao's blog, terrytao.wordpress.com, February 4, 2014.

\bibitem{Tao2}
T. Tao.
\newblock {\em Global regularity for a logarithmically supercritical hyperdissipative Navier-Stokes equation}.
\newblock Analysis and PDE, Vol. 2, No.3, 2009.


\bibitem{Vishik1}
M. Vishik.
\newblock {\em Instability and non-uniqueness in the Cauchy problem for the Euler equations of an ideal incompressible fluid. Part I}.
\newblock arXiv:1805.09426, 2018.

\bibitem{Vishik2}
M. Vishik.
\newblock {\em Instability and non-uniqueness in the Cauchy problem for the Euler equations of an ideal incompressible fluid. Part II}.
\newblock arXiv:1805.09440, 2018.




\end{thebibliography}
\end{document}